\documentclass[a4paper,12pt]{article}
\usepackage{bigints}
\usepackage{calrsfs}% for curly mathcal font %
\usepackage{soul} %strikeout
\usepackage{graphicx}
\usepackage[english]{babel}
\usepackage{hyperref}
\usepackage{amsthm,amsmath,amssymb,color}

\usepackage{fdsymbol}
\def\spad{\smalldiamond} % introduced instead of spadesuit

\newtheorem{lemma}{Lemma}
\newtheorem{coro}{Corollary}

\newtheorem{prop}{Proposition}
\newtheorem{thm}{Theorem}
\newtheorem{rema}{Remark}
\newtheorem{defn}{Definition}
\newcommand{\dd}{\,\mathrm{d}}

\def \a   {\alpha}

\def \g   {\gamma}

\def \eps {\varepsilon}
\def\E{{\mathbb{E}}}

\def\P{{\mathbb{P}}}
\def \R {{\mathbb{R}}}

\def\F{{\cal{F}}}

\title{The limit point in the Jante's law process has an absolutely continuous distribution}
\author{Edward Crane\footnote{School of Mathematics, University of Bristol, BS8 1TH, UK} {} and Stanislav Volkov\footnote{Centre for Mathematical Sciences, Lund University, Box 118 SE-22100, Lund, Sweden}
} 
\begin {document}
\maketitle
\begin{abstract}We study a stochastic model of consensus formation, introduced in 2015 by Grinfeld, Volkov and Wade, who called it a \emph{multidimensional randomized Keynesian beauty contest}. The model was generalized by Kennerberg and Volkov, who called their generalization the \emph{Jante's law process}. We consider a version of the model where the space of possible opinions is a convex body $\mathcal{B}$ in $\mathbb{R}^d$. $N$ individuals in a population each hold a (multidimensional) opinion in $\mathcal{B}$. Repeatedly, the individual whose opinion is furthest from the centre of mass of the $N$ current opinions chooses a new opinion, sampled uniformly at random from~$\mathcal{B}$. Kennerberg and Volkov showed that the set of opinions that are not furthest from the centre of mass converges to a random limit point. We show that the distribution of the limit opinion is absolutely continuous, thus proving the conjecture made after Proposition~3.2 in Grinfeld et al.
\end{abstract}



\section{Introduction}\label{Intro} \subsection{The multidimensional randomized Keynesian beauty contest}
Let $N\ge 3$ and $d \ge 1$ be integers.  Let $\lambda$ denote Lebesgue measure on $\mathbb{R}^d$ and let $\mathcal{B}$ be a \emph{convex body}, i.e. a closed convex subset of $\mathbb{R}^d$ with non-empty interior. In particular, $ \lambda(\mathcal{B}) > 0$. For the moment, while describing the model studied in \cite{GVW}, we will also assume that $\mathcal{B}$ is bounded so that $\lambda(\mathcal{B}) < \infty$, but later we will drop this assumption. We say that a random variable is uniform on $\mathcal{B}$, i.e. it is a $U(\mathcal{B})$ random variable when its distribution is the probability measure $(1/\lambda(\mathcal{B}))\lambda|_{\mathcal{B}}$. (This notion requires $\lambda(\mathcal{B}) < \infty$, of course.)

A discrete time Markov process $X(t)_{t=0}^{\infty}$, taking values in $([0,1]^d)^N$, called the \emph{multidimensional randomized Keynesian beauty contest} was studied in \cite{GVW}. In the sequel \cite{KV1} a generalization of the process was studied which the authors called the \emph{Jante's law process}\footnote{The \emph{law of Jante} is a literary caricature of the virtue in Scandinavian culture of not standing out from the crowd. It appeared in Aksel Sandemose's satirical novel \emph{En flyktning krysser sitt spor}, (A fugitive crosses his tracks), published in 1933.} (see \S\ref{SS: related models}).  In the present paper the process $X(t)$ will take values in $\mathcal{B}^N$, so our setting is only a little more general than that of \cite{GVW}; however, we prefer to retain the shorter name for the process.

If $X$ is either a finite sequence $(x_1, \dots, x_m)$ of $m \ge 1$ points of $\mathbb{R}^d$, or a set $\{x_1, \dots , x_m\}$ of $m \ge 1$ distinct points of $\mathbb{R}^d$, then we denote by $\mu(X)$ its centre of mass, i.e.
$$\mu(X) = \frac{1}{m}(x_1 + \dots + x_m)\,.$$ We denote by $\Vert x-y\Vert = \sqrt{\sum_{i=1}^d (x_i - y_i)^2}$ the Euclidean distance between $x$ and $y$ in $\mathbb{R}^d$, and by
$$
\mathrm{B}(x,\rho)=\{y\in\R^d:\ \Vert x-y\Vert <  \rho\}
$$
the open ball of radius $\rho$ centred at $x$.

Let $X(0) = (X_1(0), \dots, X_N(0)) \in\mathcal{B}^N$ be a possibly random initial state. We will always assume that a.s.~the $N$ points of $X(0)$ are distinct. An informal description of each step of the Markov chain $X(t)$ is that amongst $N$ points present at time $t$, we choose the one that is furthest from the current centre of mass, and throw it out. After that, a new point arrives, distributed uniformly in $\mathcal{B}$, and replaces the thrown-out point, so that we have again $N$ points. 

Formally, let $\left(\zeta_t\right)_{t=1}^{\infty}$ be an i.i.d.\ sequence of  random variables uniformly distributed on $\mathcal{B}$, denoted by $U(\mathcal{B})$, independent of the initial state $X(0)$.  For each time $t=0, 1, 2,\dots$ in turn, let~$j_t$ be the index of the furthest point of $X(t)$ from $\mu(X(t))$, that is,  $j_t$ is defined by
$$
\Vert X_{j_t}(t) - \mu(X(t))\Vert =  \max \left\{ \Vert X_i(t) - \mu(X(t))\Vert:  1 \le i \le N \right\}\,;
$$  
in case of a tie, choose $j_t$ uniformly at random among the tied indices\footnote{Any other way of breaking ties would also be acceptable since a.s.~ties do not occur.}. Define $X(t+1)$ by
$$
X_i(t+1) = \begin{cases} X_i(t), & \text{if $i \in \{1, \dots, N\}\setminus\{j_t\}$},\\ \zeta_t, & \text{if $i = j_t$},
\end{cases}
$$
i.e., the point of $X(t)$ that is furthest from the centre of mass of $(X_1(t), \dots, X_N(t))$ is removed and its place is taken by the new point $\zeta_t$. Note that although $X(t+1)$ and $\mu(X(t+1))$ depend continuously on the new point $\zeta_t$, they depend discontinuously on $\zeta_{t-1}$.

Let $X'(t)$ denote the {\em core} of the configuration $X(t)$, that is the \emph{set} $\{X_i(t): i \neq j_t\}$.  
\begin{defn}\label{def:conv}
We say that $X'(t)$ converges to some point $\xi\in\R^d$ if
$$
\text{for all }\eps > 0\text{ there is }t_\eps< \infty\text{ such that }  X'(t) \subset B(\xi,\eps)\text{ for all }t > t_\eps\,.
$$ 
\end{defn}
It was shown in~\cite{GVW} that $X'(t)$ converges a.s.~to a random limit point $\xi$ in the above sense; additionally, in the case $d=1$, the distribution of $\xi$ has support $[0,1]$.   It was conjectured\footnote{The conjecture is stated immediately after \cite[Proposition 3.2]{GVW}.} that the distribution of $\xi$ is, in fact, absolutely continuous with respect to the Lebesgue measure on $\mathbb{R}^d$; see the appendix for the definition. {\bf The purpose of this paper is to prove this absolute continuity conjecture in the $d$-dimensional setting described above. }
 
A.s.~for every $t \ge 0$, $X(t)$ is a sequence of $N$ distinct points, and it will be convenient to discard from the underlying probability space the null set where this fails. Hence, from now on, we work with a version of the process where for each $t$ the points of $X(t)$ are surely distinct. In particular, the core process $X'(t)$ a.s.~takes its values in the set of subsets of $\mathcal{B}$ of cardinality $M:=N-1$. Assuming $X'(t) = S$, we can ask what is the set of possible values of $\zeta_t$ which would result in $\zeta_t \in X'(t+1)$, that is, that we keep the newly sampled point in the core while throwing out one of the original points. The answer depends only on $S$, and not on the point~$X_{j_t}(t)$. The set of possible locations which could enter the core at time $t+1$ is denoted by
\begin{equation}
\mathrm{Keep}(S; \mathcal{B})\; := \;\{ x \in \mathcal{B}\,:\, \Vert x - \mu(\{x\} \cup S)\Vert < \max(\{\Vert s - \mu(\{x\} \cup S)\Vert: s \in S\}) \}.\label{def: Keep}
\end{equation}
For later use, we extend the definition~\eqref{def: Keep} to allow $\mathcal{B}$ to be unbounded.
\begin{figure}
\begin{center}
\includegraphics[scale=1]{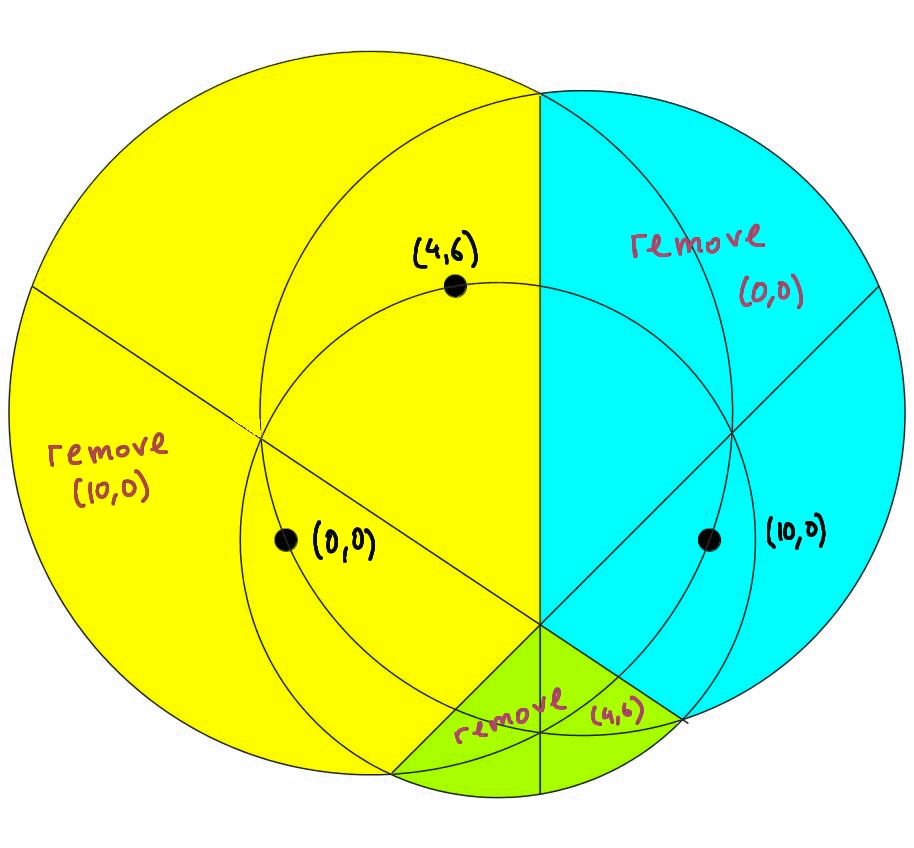}
\end{center}
\caption{$d=2$ and $M=3$. In this example, the core consists of points $(0,0)$, $(10,0)$, and $(4,6)$. A newly sampled point will join the core if and only if it lies in one of the three shaded areas; the colours indicate which point will be removed.}\label{fig1}
\end{figure}

Observe that in the one-dimensional case $\mathcal{B}$ is necessarily an interval $[a,b]\subset\R$; and then $\mathrm{Keep}(S; \mathcal{B})$ is also a non-empty subinterval of $[a,b]$. If $d\ge 2$ the geometry of this set is a little more complicated, see e.g. Figure~\ref{fig1}.
Later,  in Lemma~\ref{L: conditions to enter the core}, we will show that $\mathrm{Keep}(S; \mathcal{B})$ is in fact always the intersection of $\mathcal{B}$ with a non-empty union of at most $M$ bounded open balls which only depend on $S$. Moreover, $\mathrm{Keep}(S,\mathcal{B})$ always contains a nonempty open ball whose centre is $\mu(S)$. As a result, $\mathrm{Keep}(S;\mathcal{B})$ always has finite and positive Lebesgue measure, so it makes sense to sample a point uniformly from $\mathrm{Keep}(S;\mathcal{B})$. In fact, conditional on $X'(t) = S$ and $X'(t+1) \neq X'(t)$, the new point $\zeta_t$ is uniformly distributed on $\mathrm{Keep}(S;\mathcal{B})$.

\subsection{Time-changed models}

For convenience, we will work with a time-change of the process $X'(t)$, namely the sequence $Y(0), Y(1), Y(2), \dots$ defined by $Y(0) = X'(0)$ and for $n \ge 1$, $Y(n) = X'(t_n)$, where $t_1 < t_2 <  t_3 < \dots$ is the increasing sequence of all times $t \ge 1$ for which $X'(t) \neq X'(t-1)$.  The sequence $\left(Y(n)\right)_{n\ge 0}$ is itself a Markov chain. It is an instance of a {\bf{$\mathcal{B}$-valued Jante's law process}} (defined below). When studying $Y(\cdot)$, there is no longer any reason to insist that $\mathcal{B}$ be bounded, since we can describe the law of $Y(\cdot)$ without making any reference to $U(\mathcal{B})$ random variables. Instead, the definition involves sampling uniformly from sets of the form $\mathrm{Keep}(S;\mathcal{B})$, which always have finite positive Lebesgue measure, regardless of whether $\mathcal{B}$ is bounded. So in the following definition and for the rest of the paper {\bf we allow $\mathcal{B}$ to be unbounded}.

\begin{defn}[$\mathcal{B}$-valued Jante's law process $Y(n)_{n \in \mathbb{N}_0}$]
Let $\mathcal{B}$ be any closed convex subset of $\mathbb{R}^d$ with a non-empty interior.  Let $Y(0)$ be any (possibly random) set of $M \ge 2$ distinct points in $\mathcal{B}$. For each $n \ge 1$, sample a point $y_n$ uniformly from $\mathrm{Keep}(Y(n-1);\mathcal{B})$, conditionally independently of $(Y(i): 0 \le i < n-1)$ given $Y(n-1)$, and set
\begin{align}\label{YyR}
Y(n) := \{y_n\} \cup Y(n-1) \setminus \{r_n\}\,,
\end{align}
where $r_n$ is the (a.s.~unique) point $r \in Y(n-1)$ which maximizes $\Vert r - \mu(Y(n-1) \cup \{y_n\})\Vert$. 
\end{defn}
By definition, a $\mathcal{B}$-valued Jante's law process is a Markov chain. When working with $\mathcal{B}$-valued Jante's law processes, it is convenient to discard from our probability space the null set where for any~$n$ either there is a tie in choosing $r_n$, or $y_n \in \bigcup_{i=0}^{n-1} Y(i)$. From the definition of $\mathrm{Keep}(S;\mathcal{B})$, it is surely the case that $r_n \neq y_n$, so $Y(n) \neq Y(n-1)$.

By a slight modification of the proof in \cite{GVW} we will show that even when $\mathcal{B}$ is unbounded, the $\mathcal{B}$-valued Jante's law process $Y(\cdot)$ a.s.~converges to a random limit point $\xi$, in the sense that
$$(\forall \eps > 0)(\exists n_\eps< \infty)(n > n_\eps \implies Y(n) \subset B(\xi,\eps))\,.$$ 

Our strategy for proving that the limit point $\xi$ is an absolutely continuous random variable is to compare~$Y(\cdot)$ with another process, denoted by $Z(n)_{n \in \mathbb{N}_0}$, which is just a $\mathcal{B}$-valued Jante's law process with $\mathcal{B}=\R^d$; we call the Markov chain $Z(\cdot)$ the \emph{scale-free core process}. It eliminates the complicated effect of the boundary of $\mathcal{B}$ and has the great advantage of being invariant in law with respect to scaling and translation.\footnote{An analogous modification of the original process $X(t)$ was already used as a tool in \cite[\S3.4]{GVW}.}  It will be useful to have notation for the points added and removed at each step of the scale-free core process $Z(n)$:
$$
Z(n) = \{z_n\} \cup (Z(n-1) \setminus \{s_n\})\,,
$$ 
Also let $z_\infty$ denote the a.s.~limit of the process $Z(n)$, in the sense of Definition~\ref{def:conv}.

\subsection{Outline and results} 
In \S\ref{S: preliminaries} we collect a number of basic results, including the fact that every $\mathcal{B}$-valued Jante process almost surely converges to a random limit (Proposition~\ref{prop:conv}). We give a standalone proof because for the case where $\mathcal{B}$ is unbounded, this is not a special case of the results of \cite{GVW,KV1}, and because the required lemmas are needed again later in the paper.

In \S\ref{S: uniform geometry} we discuss a property of convex bodies called \emph{uniform geometry}. This property is used later in the paper to avoid some geometric complications. Bounded convex sets have uniform geometry. We reduce the absolute continuity of the distribution of the limit point of a $\mathcal{B}$-valued Jante process to the case where $\mathcal{B}$ is bounded.

Our first continuity result is for the scale-free core process:
\begin{thm}\label{thm: scale-free limit is continuous}
For an arbitrary deterministic initial condition $Z(0) = S$, where $S$ is any set of~$M$ distinct points in $\mathbb{R}^d$, the distribution $\pi_S$ of $z_\infty $ is absolutely continuous. 
\end{thm}

Our proof of Theorem~\ref{thm: scale-free limit is continuous} is less straightforward than one might expect. The main difficulty is to show that a.s.~all of the original points are eventually removed; we isolate this statement as Proposition~\ref{prop: exodus} and prove it in \S\ref{S: exodus} using a supermartingale argument, a trick using translational symmetry, and a compactness argument. Then in \S\ref{S: continuous limit for scale-free process} we complete the proof of Theorem~\ref{thm: scale-free limit is continuous} by showing that the distribution of the limit is a mixture of absolutely continuous distributions. Each distribution in the mixture is seen to be absolutely continuous because it is a convolution with a uniform distribution over a set of positive Lebesgue measure.

In \S\ref{S: coupling succeeds} we use another supermartingale argument and a coupling argument to deduce the absolute continuity of the distribution of the limit point for any $\mathcal{B}$-valued Jante's law process. 
\begin{thm}\label{T: original}
Let $\mathcal{B} \subset \mathbb{R}^d$ be any convex body, and let $S$ be any set of  $M$ distinct points of $\mathcal{B}$. Let $\xi$ be the limit point of the $\mathcal{B}$-valued Jante's law process $Y(\cdot)$ started at $Y(0) = S$. Then the distribution of $\xi$ is absolutely continuous. 
\end{thm}

Of course this also shows that the distribution of the random limit $\xi$ of the original core process $X'(t)$ started at $X'(0) = S$ has an absolutely continuous distribution. Theorem~\ref{T: original} covers the original setting of~\cite{GVW} where $\mathcal{B} = [0,1]^d$, so it resolves the absolutely continuous distribution conjecture made in that paper.  We remark that if $X'(0)$ is random, then the distribution of the limit point $\xi$ is still absolutely continuous since it is a mixture of such distributions. We also remark that $\mathbb{P}(\xi \in \partial\mathcal{B}) = 0$ because $\lambda(\partial\mathcal{B}) = 0$ for any convex body $\mathcal{B} \subseteq \mathbb{R}^d$.

In Appendix 3 we discuss the main obstacle to extending the Theorem~\ref{T: original} to a larger class of subsets $\mathcal{B} \subset \mathbb{R}^d$. For example, if $\mathcal{B}$ is a closed and bounded subset of $\mathbb{R}^d$ whose boundary is a smooth compact hypersurface embedded in $\mathbb{R}^d$, then we can extend the definition of the $\mathcal{B}$-valued Jante process and check using the results of \cite{KV1} that it converges almost surely to a random limit $y_\infty \in \mathcal{B}$. But we do not have a proof that $\mathbb{P}(y_\infty \in \partial \mathcal{B}) = 0$, which is certainly necessary for the absolute continuity of the distribution of $y_\infty$. On the other hand, we show in Lemma~\ref{L: non-convex continuity} that if $\mathbb{P}(y_\infty \in \partial \mathcal{B}) = 0$ then the distribution of $y_\infty$ is absolutely continuous. 

\subsection{Related models}\label{SS: related models} In the paper \cite{KV1}, the authors generalized the original model of \cite{GVW} by allowing the common distribution of the new points $\zeta_t$ to be an arbitrary absolutely continuous distribution on $\mathbb{R}^d$. They also allowed for more than one point to be replaced simultaneously. It was shown that if the distribution of~$\zeta_t$ has bounded support, together with some extra regularity assumptions, then a.s.~a limit point~$\xi$ exists. However, although the support of~$\xi$ is certainly a subset of the support of~$\zeta$, it may not be the whole of the support of~$\zeta$. For example, if~$\zeta_t$ is uniform on $[0,1] \cup [2,3]$ and a.s.~all coordinates of $X(0)$ lie in $[0,1]$, then $X'(t) \subset [0,1]$ for all $t$, and after passing to the subsequence of times at which $\zeta_t \in [0,1]$, we recover the original process with $\zeta \sim U([0,1])$. In particular, in this case, the support of $\xi$ is $[0,1]$. In case $\zeta$ has unbounded support, \cite{KV1} shows that a.s.~either $X'(t)$ converges to a random limit point $\xi$ \emph{or} $X'(t) \to \infty$, i.e. $\inf\{|x|: x \in X'(t)\} \to \infty$.  However, no examples of distributions of $\zeta_t$ are known to satisfy $X'(t) \to \infty$ with positive probability. For an arbitrary common distribution of the random variables $\zeta_t$, and an arbitrary distribution of $X(0)$ on $\mathbb{R}^d$, such that there is a positive probability that $X(t)$ converges, one can ask whether the distribution of the limit $\xi$ conditional on its existence is absolutely continuous with respect to the distribution of the $\zeta_t$. We do not address this more general continuity problem in this paper, although it might be possible to extend our techniques to other cases where $\zeta_t$ has an absolutely continuous distribution with a sufficiently well-behaved density.

In another paper~\cite{KV2},  the authors studied asymptotic properties of a modified version of Jante's law process. In this model, called the \emph{$p$-contest}, at each moment of time  the point farthest from $p\mu$, where $\mu$ is the current centre of mass, is removed and replaced by an independent
$\zeta$-distributed point; note that the case $p=1$ would correspond to the original Jante's law process. Finally, in~\cite{KV3}, a local version of the process was considered. In this version, $N$ vertices are placed on a circle, so that each vertex has exactly two neighbours. To each vertex assign a real number, called {\em fitness}. At each unit of time, the vertex whose fitness deviates most from the average of the fitnesses of its two immediate neighbours has its fitness replaced by a random value drawn independently according to some distribution $\zeta$. The authors showed that if $\zeta$ has a uniform or a discrete uniform distribution, all fitnesses but at most one converge to the same limit. We would also like to mention that there is a related continuous-time model, called ``Brownian bees'', which involves a fixed number $N$ of particles which perform independent Brownian motions, and branch from time to time simultaneously with the removal of the particle most distant from the origin, thus keeping the total number of particles constant; see e.g.~\cite{ABL} and \cite{SSM}. It would be in the spirit of the current paper to identify extrema in the Brownian bees from the centre of mass of the swarm, rather than from the origin; and have exponential killing of the extreme bee with introducing a new bee at a random location.

\section{Preliminaries}\label{S: preliminaries}

For any set $X \subset \mathbb{R}^d$ we denote by $\mathrm{Conv}(X)$ the convex hull of $X$, i.e.\ the smallest convex set containing $X$. We denote by $\partial X$ the boundary of $X$ and by $X^\circ$ the interior of $X$. For any $\delta > 0$ we denote by $N_\delta(X)$ the $\delta$-neighbourhood of $X$, i.e.\ the set of points whose Euclidean distance to $X$ is less than $\delta$. 

Some basic facts about continuity of $\mathbb{R}^d$-valued random variables and continuity of Borel probability measures on $\mathbb{R}^d$ are collected in the appendix. The main two facts that we will use are that any mixture of absolutely continuous random variables is absolutely continuous and that if an absolutely continuous random variable is conditioned on an event of positive probability then the conditioned random variable is absolutely continuous. We will see a simple argument that uses these two facts in the proof of Lemma~\ref{L: reduction to bounded case} at the end of \S\ref{S: uniform geometry}.

For any set $X = \{x_1, \dots, x_M\} \subset \mathbb{R}^d$ of cardinality $M = N-1$ we let
$$ \Sigma(X)  = M \mu(X) = \sum_{i=1}^M x_i$$ 
and we define the following geometric functionals of $X$: 
\begin{align}\label{dDA}
\begin{array}{rclrcl}
F(X)&=& \displaystyle \sum_{1 \le i < j \le M} \Vert x_i - x_j\Vert^2\,=\,M\sum_{i=1}^M \Vert x_i - \mu(X)\Vert^2\,,
&d(X) &=& \displaystyle \min_{i \neq j} \Vert x_i - x_j\Vert\,, \\
A(X) &=& \max_{w \in X} \Vert w - \mu(X)\Vert\,,
&D(X) &=& \max_{i \neq j} \Vert x_i - x_j\Vert\,.
\end{array}
\end{align}
We call $F(X)$ the {\em moment of inertia} of $X$. The functional $F$ serves as a Lyapunov function for the Jante's law processes $Y(\cdot)$ and $Z(\cdot)$. 
One can easily verify the following inequalities among the geometric functionals:
\begin{align}\label{useful}
\begin{split}
 \frac{M}{M-1}\, A(X)  &\le D(X) \le 2A(X),\\
 \frac{M(M-1)}{2}\,d(X)^2 &\le F(X) \le \frac{M(M-1)}{2}\, D(X)^2, \\
 \frac{M^2}{(M-1)}\,A(X)^2 &\le F(X) \le M^2A(X)^2.
\end{split}
\end{align}
 
\begin{lemma}\label{L: about F}
Let $X = \{x_1, \dots, x_M\}$ be any set of $M \ge 2$ distinct points in $\mathbb{R}^d$, and suppose $z \in \mathbb{R}^d\setminus X$. Then \begin{eqnarray}\label{E: decrease of F} F(X) - F\left(\{z\} \cup X \setminus \{x_j\}\right) & = & (M+1)\left(\left\Vert x_j - \frac{z+\Sigma(X)}{M+1}\right\Vert^2 - \left\Vert z - \frac{z+\Sigma(X)}{M+1}\right\Vert^2\right)\\ \label{E: closer to the mean of the others} & = & (M-1)\left(\left\Vert x_j -  \frac{\Sigma(X) - x_j}{M-1}\right\Vert^2 - \left\Vert z - \frac{\Sigma(X) - x_j}{M-1} \right\Vert^2\right)\,.\end{eqnarray}
From~\eqref{E: decrease of F} it follows that a point  $z \in \mathcal{B}\setminus X$ belongs to $\mathrm{Keep}(X;\mathcal{B})$ if and only if there exists $x \in X$ such that $F(\{z\} \cup (X \setminus \{x\})) < F(X)$. Moreover, if $z \in \mathrm{Keep}(X;\mathcal{B})$ then any choice of $j$ which maximizes $\Vert x_j - \frac{\Sigma(X) + z}{M+1}\Vert$ is also a choice which minimizes
$F(\{z\} \cup (X \setminus\{x_j\}))$. 
\end{lemma} 
\begin{proof}
We have
\begin{eqnarray*} 
F(X) - F(\{z\} \cup (X \setminus \{x_j\})) & = & \sum_{i: i \neq j} \Vert x_j - x_i \Vert^2 - \sum_{i: i \neq j} \Vert z - x_i\Vert^2 \\
 & = & \sum_{i: i \neq j} (\Vert x_j \Vert^2 - \Vert z \Vert^2 - 2(x_j - z)\cdot x_i) \\ & = & (M-1)(x_j - z)\cdot(x_j + z) - 2(x_j- z)\cdot(\Sigma(X) - x_j)\\ 
 & = & (x_j - z)\cdot( (M+1)(x_j + z) -2z - 2\Sigma(X))\,.
   \end{eqnarray*}  
On the other hand, expanding the squared distances on the RHS of~\eqref{E: decrease of F} as dot products, the terms $\pm (M+1)\left\Vert\frac{z + \Sigma(X)}{M+1}\right\Vert^2$ cancel, and we are left with the same quantity:
\begin{equation*}(M+1) ( \Vert x_j\Vert^2 - \Vert z\Vert^2) - 2 (x_j-z)\cdot(z + \Sigma(X))  =  (x_j - z)\cdot((M+1)(x_j+z) - 2z - 2\Sigma(X))\,. \end{equation*} The proof of~\eqref{E: closer to the mean of the others} is a similar calculation.
The remaining statements follow immediately from~\eqref{E: decrease of F}.
\end{proof}
\begin{coro}[Essentially Lemma~2.1 in~\cite{GVW}]\label{cor: F decreases}
Let $Y(\cdot)$ be a $\mathcal{B}$-valued Jante's law process with $M \ge 2$ (where $\mathcal{B}$ may be~$\R^d$). Then $F(Y(n)) < F(Y(n-1))$ for each $n\ge 1$.
\end{coro}

\begin{lemma}\label{L: conditions to enter the core}
Let $X = \{x_1, \dots, x_M\} \subset \mathcal{B}$, where $M \ge 2$.   Then 
\begin{equation}\label{E: Keep as a union of balls}
\mathrm{Keep}(X;\mathcal{B}) = \mathcal{B} \cap \bigcup_{j=1}^M B\left(\frac{\Sigma(X)-x_j}{M-1}\,,\, \frac{M}{M-1}\Vert x_j - \mu(X)\Vert\right)\,,
\end{equation}
\end{lemma}

\begin{proof}
Let $X = \{x_1, \dots, x_M\}$. From the definition~\eqref{def: Keep}, for any point $z \in \mathbb{R}^d$, $z \in \mathrm{Keep}(X;\mathbb{R}^d)$ if and only if there exists $j \in \{1, \dots, M\}$ such that
\begin{equation} 
\label{E: condition to enter core}
\left\Vert z - \frac{z + \Sigma(X)}{M+1}\right\Vert < \left\Vert  x_j - \frac{z + \Sigma(X)}{M+1}\right\Vert\,.
\end{equation}
Straightforward algebraic manipulation shows that each of the following three inequalities is equivalent to~\eqref{E: condition to enter core}.
\begin{align}\label{E: simple}
\begin{split}
\left\Vert M z -  \Sigma(X) \right\Vert &< \left\Vert(M+1) x_j - z - \Sigma(X)\right\Vert
\\
%\end{equation} 
%\begin{equation*}\label{E: F decreases}
\sum_{i: i \neq j} \Vert z - x_i\Vert^2 &< \sum_{i: i \neq j} \Vert x_j - x_i\Vert^2
\\
%\end{equation*} 
%\begin{equation}
%\label{E: Stewart} 
\left\Vert z - \frac{\Sigma(X) - x_j}{M-1}\right\Vert &< \left\Vert x_j - \frac{\Sigma(X) - x_j}{M-1}\right\Vert = \frac{M}{M-1}\Vert x_j - \mu(X)\Vert\,.
\end{split}
\end{align}
The expression~\eqref{E: Keep as a union of balls} follows from the last of these inequalities. 
\end{proof}
\begin{lemma}\label{L: ball containing Keep}
Let $X = \{x_1, \dots, x_M\} \subset \mathcal{B}$, where $M \ge 2$.   Then 
\begin{equation}\label{E: sandwich1} 
 \mathrm{Keep}(X;\mathcal{B}) \subseteq \mathcal{B} \cap B\left(\mu(X),\frac{M+1}{M-1}A(X)\right)\,.
\end{equation}
\end{lemma}
\begin{proof}
 For every $j \in \{1, \dots, M\}$ we have \begin{equation}\label{E: Fred} \left\Vert \mu(X) - \frac{\Sigma(X) - x_j}{M-1}\right\Vert = \frac{1}{M-1}\Vert x_j - \mu(X)\Vert \le \frac{1}{M-1}A(X)\,.\end{equation}
Suppose $z \in \mathrm{Keep}\left(X;\mathbb{R}^d\right)$ and let $j$ be such that the equivalent inequalities~\eqref{E: condition to enter core}-\eqref{E: simple}
%\eqref{E: Stewart} 
 hold. Then by~\eqref{E: simple} 
%\eqref{E: Stewart} 
 we have
\begin{equation}\label{E: Ginger} \left\Vert z - \frac{\Sigma(X) - x_j}{M-1}\right\Vert  < \frac{M}{M-1}\Vert x_j - \mu(X)\Vert \le \frac{M}{M-1}A(X)\,.\end{equation}
Summing the inequalities~\eqref{E: Fred}~and~\eqref{E: Ginger}, and using the triangle inequality,  we obtain
$$ 
\Vert z - \mu(X)\Vert  \le \left\Vert z - \frac{\Sigma(X) - x_j}{M-1}\right\Vert  +  \left\Vert \mu(X) - \frac{\Sigma(X) - x_j}{M-1}\right\Vert < \frac{M+1}{M-1}\cdot A(X)\,,
$$
i.e. $z \in B\left(\mu(X),\frac{M+1}{M-1}\, A(X)\right)$. 
\end{proof}

\begin{lemma}\label{L: ball inside Keep}
Let $X = \{x_1, \dots, x_M\} \subset \mathcal{B}$, where $M \ge 2$ and the $x_i$ are distinct.   Then 
\begin{equation}\label{E: sandwich2} 
\mathcal{B} \cap B\left(\mu(X), A(X)\right) \subseteq \mathrm{Keep}(X;\mathcal{B})\,. 
\end{equation}
\end{lemma}
\begin{proof}
Let $x_j$ be a point of $X$ furthest from $\mu(X)$. Then $\mu(X)$ lies on the line segment joining $\frac{\Sigma(X) - x_j}{M-1}$ and $x_j$, since $\mu(X)=\left(1-\frac1M\right)\left(\frac{\Sigma(X)-x_j}{M-1}\right) + \frac{1}{M} x_j$. We have
 $\Vert \mu(X) - x_j \Vert = A(X)$ and $\left\Vert \frac{\Sigma(X) - x_j}{M-1} - x_j\right\Vert  = \left\Vert \frac{M\mu(X)}{M-1} - \frac{M}{M-1}x_j\right\Vert  = \frac{M}{M-1}A(X)$. Therefore $$B(\mu(X),A(X)) \subseteq B\left(\frac{\Sigma(X) - x_j}{M-1}, \frac{M}{M-1}\Vert x_j - \mu(X)\Vert\right)\,.
 $$
Indeed, the boundary spheres of these two balls are tangent at the point $x_j$. The result then follows from Lemma~\ref{L: conditions to enter the core}. 
\end{proof}
\begin{lemma}\label{L: other balls inside Keep}
Let $X = \{x_1, \dots, x_M\} \subset \mathcal{B}$, where $M \ge 2$ and the $x_i$ are distinct. For each $x_i \in X$ we have
$$ \mathcal{B} \cap B\left(\mu(X \setminus \{x_i\}), \frac{M}{M+1}A(X)\right) \subseteq \mathrm{Keep}(X; \mathcal{B})\,. $$
\end{lemma}
\begin{proof}
 $\mu(X \setminus \{x_i\}) = \frac{\Sigma(X) - x_i}{M-1}$ so $\Vert\mu(X) - \mu(X \setminus\{x_i\})\Vert = \frac{1}{M-1}\Vert \mu(X) - x_i \Vert$.  Therefore by Lemma~\ref{L: ball inside Keep} we have
 $$ \mathcal{B} \cap B\left(\mu(X \setminus\{x_i\}), A(X) - \frac{1}{M-1}\Vert \mu(X) - x_i \Vert \right) \subseteq \mathrm{Keep}(X; \mathcal{B})\,.$$
 On the other hand from Lemma~\ref{L: conditions to enter the core} we have
 $$ \mathcal{B} \cap B\left(\mu(X \setminus\{x_i\}), \frac{M}{M-1}\Vert \mu(X) - x_i \Vert \right) \subseteq \mathrm{Keep}(X; \mathcal{B})\,.$$
 The maximum of the two radii defining these sets is at least $\frac{M}{M+1}A(X)$. 
\end{proof}

\begin{lemma}\label{L: good decrease in F}
Let $X = \{x_1, \dots, x_M\} \subset \mathcal{B}$, where $M \ge 2$.  Suppose $z \in B\left(\mu(X),\alpha A(X)\right)$ for some $0 < \alpha < 1$, and $z \not\in X$. Let $x_k$ be a point of $X$ farthest from $(z+\Sigma(X))/(M+1)$. Then 
\begin{equation}\label{E: definite decrease of F}
F(X) - F(\{z\} \cup (X \setminus\{x_k\})) \ge (1-\alpha)((M-1)\alpha + M+1)A(X)^2\, > \frac{1-\alpha}{M}\, F(X)\,.
\end{equation}
\end{lemma}

\begin{proof}
Let $x_j$ be a point of $X$ farthest from $\mu(X)$, so that $\Vert x_j - \mu(X)\Vert = A(X)$ and 
$$\left\Vert x_j - \frac{\Sigma(X) - x_j}{M-1}\right\Vert = \frac{M}{M-1}A(X)\,.$$
By the triangle inequality,
$$\left\Vert z - \frac{\Sigma(X) - x_j}{M-1} \right\Vert \le \Vert z - \mu(X) \Vert + \left\Vert \mu(X) - \frac{\Sigma(X) - x_j}{M-1}\right\Vert \le \alpha A(X) + \frac{1}{M-1}A(X)\,.$$
By Lemma~\ref{L: about F} and equation~\eqref{E: closer to the mean of the others} we have
\begin{eqnarray*}F(\{z\} \cup (X \setminus \{x_k\})) & \le & F(\{z\} \cup (X \setminus \{x_j\})) \\ & = & F(X) - (M-1) \left(\left\Vert x_j - \frac{\Sigma(X) - x_j}{M-1}\right\Vert^2 - \left\Vert z - \frac{\Sigma(X) - x_j}{M-1} \right\Vert^2\right)\\
& \le & F(X) - (M-1)\left(\left(\frac{M}{M-1}\right)^2 - \left(\alpha + \frac{1}{M-1}\right)^2  \right)A(X)^2 \\ 
 & = & F(X) - (1-\alpha)((M-1)\alpha + M+1))A(X)^2\,.
\end{eqnarray*}
The final strict inequality in~\eqref{E: definite decrease of F} follows from the final inequality in~\eqref{useful}.\footnote{In fact one  could obtain $F(X) - F(\{z\} \cup (X \setminus \{x_k\})) \ge \frac{(1-\alpha)(M+1)}{M^2} F(X)$, but the simpler lower bound given in~\eqref{E: definite decrease of F} will suffice for the applications.}

\end{proof}

The following estimate is well-known (see e.g. \cite[Lemma 2.6]{Leonardi}).
\begin{lemma}\label{L: volume comparison}
Let $\mathcal{B} \subset \mathbb{R}^d$ be a convex body. Let $x \in \mathcal{B}$ and $0 < r_1 < r_2$. Then 
$$\frac{\lambda(\mathcal{B} \cap B(x,r_1))}{\lambda(\mathcal{B} \cap B(x,r_2))} \ge \left(\frac{r_1}{r_2}\right)^d\,.$$
\end{lemma}

\begin{proof}
Since $d$-dimensional Lebesgue measure $\lambda$ is translation-invariant, we may assume w.l.o.g.~that $x=0$. Since $\mathcal{B}$ is convex and contains $0$, the dilation $(r_2/r_1)\mathcal{B}$ contains $\mathcal{B}$, so
$$ \mathcal{B} \cap B(0,r_2) \subseteq \left(\frac{r_2}{r_1}\right)(\mathcal{B} \cap B(0,r_1))\,.$$ Comparing volumes yields the desired inequality.
\end{proof}
\begin{coro}[Adapted from Lemma~2.3 and Lemma~2.5 in~\cite{GVW}]\label{cor: significant decrease of F}
Let $Y(0), Y(1), Y(2), \dots$ be a $\mathcal{B}$-valued Jante's law process, where $\mathcal{B}$ is any convex subset of $\mathbb{R}^d$ with non-empty interior. Let~$\F_n$ be the $\sigma$-algebra generated by $Y(i)_{i \le n}$. Then for every $n \ge 1$, \begin{equation}\label{E: probability of definite decrease in F}
    \P\left(F(Y(n+1))- F(Y(n)) < -\frac{1}{4M} F(Y(n))\,\mid\, \mathcal{F}_n\right) \,\ge\, 4^{-d}\,.\end{equation}
Consequently,
\begin{align}
\label{E: expectation of decrease in F}
\mathbb{E}(\log F(Y(n+1)) - \log F(Y(n)) \,|\, \F_n) \le -\frac{1}{4M}\,4^{-d}\,.
\end{align}
\end{coro}

\begin{proof}
Suppose $Y(n) = X$. Since $\mathcal{B}$ is convex, $\mu(X) \in \mathcal{B}$. 
The new point $y_{n+1}$ in $Y(n+1)$ is distributed uniformly in $\mathrm{Keep}(X; \mathcal{B})$, and by Lemmas~\ref{L: ball containing Keep} and~\ref{L: ball inside Keep} we have
$$\mathcal{B} \cap B\left(\mu(X), A(X)\right) \subseteq \mathrm{Keep}(X; \mathcal{B}) \subseteq \mathcal{B} \cap B\left(\mu(X), \frac{M+1}{M-1}A(X)\right)\,.$$  
Since $M \ge 2$, we have $\frac{M+1}{4(M-1)}A(X) \le A(X)$, so $$\mathbb{P}\left(y_{n+1} \in \mathcal{B} \cap B\left(\mu(X),\frac{M+1}{4(M-1)}A(X)\right)\right) \ge \frac{\lambda\left(\mathcal{B} \cap B\left(\mu(X),\frac{M+1}{4(M-1)}A(X)\right)\right)}{\lambda\left(\mathcal{B} \cap B\left(\mu(X), \frac{M+1}{M-1}A(X)\right)\right)}\,\ge 4^{-d}\,,$$ where the second inequality follows from Lemma~\ref{L: volume comparison}. When $y_{n+1} \in B\left(\mu(X),\frac{M+1}{4(M-1)}A(X)\right)$, we may take $X = Y(n)$ and $\alpha = \frac{M+1}{4(M-1)}$ in Lemma~\ref{L: good decrease in F}, so that inequality~\eqref{E: definite decrease of F} becomes
\begin{align}\label{eq: FF12}
F(Y(n)) - F(Y(n+1)) \;>\; \frac{3M-5}{4M(M-1)}\,F(Y(n)) \;\ge\; \frac{1}{4M}\, F(Y(n))\,.
\end{align}
% 3M-5 \ge M-1  iff M \ge 2
The final consequence follows from the fact that $F(Y(n+1))/F(Y(n)) \le 1$ (by Corollary~\ref{cor: F decreases}) and the inequality $\log x \le x-1$ for $0 < x \le 1$.
\end{proof}

Inequality~\eqref{E: probability of definite decrease in F}  implies that a.s.~$Y(n)$ converges exponentially fast.
\begin{lemma}[essentially Lemma~2.5 in \cite{GVW}]\label{L: exponential shrinking}
Let $N\ge 3$. For any fixed choice of $Y(0)$, there exists $\a > 0$ such that, a.s., $D(Y(n)) \le  e^{-\a n}$ for all sufficiently large $n$.
\end{lemma}
\begin{proof}
 From Corollary~\ref{cor: F decreases} and~\eqref{E: probability of definite decrease in F} it follows that $\mathbb{E}[F(Y(n))| F(Y(0))] \le \gamma^n F(Y(0))$, where $\gamma = 1-\frac{4^{1-d}}{M} < 1$. 
 Hence $\mathbb{P}(F(Y(n)) \ge \gamma^{n/2} F(Y(0)) ) \le \gamma^{n/2}$, and therefore a.s.\ for all sufficiently large $n$ we have $F(Y(n)) < \gamma^{n/2} F(Y(0))$. The lemma follows using the inequalities~\eqref{useful} to bound $D(Y(n))$ above in terms of $\sqrt{F(Y(n))}$; we may take any $\alpha > 0$ such that $e^{-\alpha} < \gamma^{1/4}$.   
\end{proof}
From the preliminary results established above it is now easy to recover the almost sure convergence result \cite[Theorem 1.1]{GVW} (for $\mathcal{B} = [0,1]^d$) extended to the general setting of a $\mathcal{B}$-valued Jante's law process $Y(n)$; we state it below as Proposition~\ref{prop:conv}. For the case of a bounded convex body $\mathcal{B}$, this is in fact a special case of \cite[Theorem 2]{KV1}: the regularity assumption required by that theorem is satisfied by the uniform distribution on $\mathcal{B}$, as an immediate consequence of Lemma~\ref{L: volume comparison}. In the case of unbounded $\mathcal{B}$ the exponential decay of $D(Y(n))$ ensures that the process will not escape to infinity, a fact which was also implicitly used in~\cite[\S3.4]{GVW}.
\begin{prop}[essentially Theorem~1.1 in \cite{GVW}]
\label{prop:conv}
Let $d\ge 1$ and $M=N-1\ge 2$. Let $Y(0)$ consist of~$M$ distinct points in $\mathcal{B}$. Then there exists a random $\xi\in \mathcal{B}$ such that $Y(n)$ converges to $\xi$ a.s.\ in the sense of Definition~\ref{def:conv}.
\end{prop}

It will be useful to know that if $F(Y(n))$ is small, then with high probability the $\mathcal{B}$-valued Jante's law process never travels too far from $\mu(Y(n))$. 
\begin{lemma}\label{L: relative tightness of xi}
Fix any $\epsilon > 0$ and $n\in{\mathbb N}$. Let $\gamma = 1 - 4^{1-d}/M$ and  $n_0(\epsilon) = \left\lceil 2 \log_\gamma(\epsilon (1 - \sqrt{\gamma}))\right\rceil$. Denote
$$
{\cal E}_j=\left\{\Vert \mu(Y(n+j)) - \mu(Y(n)) \Vert \le \frac{2}{M\sqrt{M-1}}\left(n_0(\epsilon) + \frac{1}{1-\gamma^{1/4}}\right)\sqrt{F(Y(n))}\right\}
$$
Then
$$
\P\left(\bigcap_{j=1}^\infty {\cal E}_j \mid Y(n)\right)\ge 1-\epsilon\,.
$$
As a result,
$$
\P\left(\left.\Vert \xi - \mu(Y(n)) \Vert \le \frac{2}{M\sqrt{M-1}}\left(n_0(\epsilon) + \frac{1}{1-\gamma^{1/4}}\right)\sqrt{F(Y(n))} \,\right|\, Y(n)\right)\ge 1-\epsilon
$$
as well.
%Then for any $n \in \mathbb{N}$, conditional on $Y(n)$, with probability at least $1-\epsilon$ we have for all $j \ge 1$ that  $$ \Vert \mu(Y(n+j)) - \mu(Y(n)) \Vert \le \frac{2}{M\sqrt{M-1}}\left(n_0(\epsilon) +\frac{1}{1-c^{1/4}}\right)\sqrt{F(Y(n))}\,, $$ and hence $$ \Vert \xi - \mu(Y(n)) \Vert \le \frac{2}{M\sqrt{M-1}}\left(n_0(\epsilon) +\frac{1}{1-c^{1/4}}\right)\sqrt{F(Y(n))}\,, $$
\end{lemma}

\begin{proof}
From Lemma~\ref{L: ball containing Keep} and \eqref{useful}, we have for every $t \ge 1$
\begin{eqnarray*}
\Vert \mu(Y(t)) - \mu(Y(t-1)) \Vert & = & \left\Vert \frac{y_t - r_t}{M}\right\Vert \\ & \le & \left(1 +\frac{M+1}{M-1}\right)\frac{A(Y(t-1))}{M} = \frac{2}{M-1}A(Y(t-1))\\ & \le & \frac{2}{M\sqrt{M-1}}\sqrt{F(Y(t-1))}\,. \end{eqnarray*}
Since $F(Y(n+i)) \le F(Y(n))$ for each $i = 0, \dots, n_0(\epsilon) - 1$, by applying the triangle inequality to the cases $t=n+1, \dots, n+n_0(\epsilon)$ of the above displayed inequalities, we obtain
$$
\Vert \mu(Y(n+ n_0(\epsilon))) - \mu(Y(n))\Vert \le n_0(\epsilon) \frac{2}{M\sqrt{M-1}}\sqrt{F(Y(n))}\,.
$$
Now note that $n_0(\epsilon)$ is the least natural number $k$ such that $\sum_{i=k}^\infty \gamma^{i/2} \le \epsilon$. As we showed in the proof of Lemma~\ref{L: exponential shrinking}, $\mathbb{P}(\,F(Y(n+i) > \gamma^{i/2}F(Y(n)))\,) \le \gamma^{i/2}$, so with probability at least $1-\epsilon$ it holds for all $i \ge n_0(\epsilon)$ that $F(Y(n+i)) \le \gamma^{i/2} F(Y(n))$. Suppose that this event occurs. Applying the triangle inequality to the cases $t \ge n+n_0(\epsilon)+1$, we may bound the relevant increments by a geometric series. For every $j > n_0(\epsilon)$ we have
$$ 
\Vert \mu(Y(n + j)) - \mu(Y(n + n_0(\epsilon)) \Vert 
\le \left(\sum_{i = n_0(\epsilon)}^{j} \gamma^{i/4}\right) \frac{2}{M\sqrt{M-1}} \sqrt{F(Y(n))}
$$ 
For all $j > n_0(\epsilon)$ the geometric series is bounded by $1/(1-\gamma^{1/4})$. Applying the triangle inequality yet again and using $\xi = \lim_{j \to \infty} Y(n+j)$ gives the claimed bounds.  
\end{proof}

\section{Reduction to the case of uniform geometry}\label{S: uniform geometry}
Following the terminology of \cite[Ch. 3]{Leonardi}, a convex body $\mathcal{B}$ is said to have \emph{uniform geometry} when there exists some $r_0 > 0$ such that
\begin{equation}\label{E: uniform geometry} 
b(r_0) := \inf_{x \in \mathcal{B}} \lambda(B(x,r_0) \cap \mathcal{B}) > 0
\end{equation}
(see \cite[Ch. 3]{Leonardi}.)  Since the volume $\lambda(B(x,r_0) \cap \mathcal{B})$ depends continuously on $x$, when $\mathcal{B}$ is compact the infimum is achieved and is positive because $\mathcal{B}$ is the closure of its interior. That is, every bounded convex body has uniform geometry. Examples of unbounded convex bodies with bounded geometry are  $\mathbb{R}^d$ itself, and any convex body which is the intersection of finitely many closed half-spaces. However, not every unbounded convex body has uniform geometry, see for example \cite[Example 3.12]{Leonardi}. 

Suppose $\mathcal{B}$ is a convex body of uniform geometry, with $r_0$ and $b(r_0) > 0$ as in~\eqref{E: uniform geometry}. Let $V(d)$ denote the volume of the unit ball in $\mathbb{R}^d$, and define 
\begin{equation} 
\label{E: c} 
c = c(\mathcal{B}, r_0) := \frac{b(r_0)}{V(d) r_0^d}\,. \end{equation}  
Note that in the case $\mathcal{B} = \mathbb{R}^d$ we have $c = 1$, and in all other cases we have $c \le 1/2$. For the case $\mathcal{B} = [0,1]^d$ which was studied in \cite{GVW}, we have $c = 2^{-d}$. By Lemma~\ref{L: volume comparison}, using 
\eqref{E: uniform geometry}  and \eqref{E: c},
for all $r \in (0,r_0)$ and all $x \in \mathcal{B}$ we have
\begin{equation}\label{E: uniform geometry 2} \frac{\lambda(\mathcal{B} \cap B(x,r))}{\lambda(B(x,r))} \ge c > 0\,. \end{equation}

Recalling Lemma~\ref{L: ball containing Keep} and Lemma~\ref{L: ball inside Keep}, for $X = \{x_1, \dots, x_M\} \subset \mathcal{B}$ with $M \ge 2$ we have
\begin{equation}\label{E: two inclusions}\mathcal{B} \cap B\left(\mu(X), A(X)\right) \subseteq \mathrm{Keep}(X; \mathcal{B}) \subseteq \mathcal{B} \cap B\left(\mu(X), \frac{M+1}{M-1} A(X)\right).\end{equation}
It follows from~\eqref{E: uniform geometry 2} and ~\eqref{E: two inclusions} that when $A(X) \le r_0$, $$ \frac{\lambda(\mathrm{Keep}(X;\mathcal{B}))}{\lambda\left(B\left(\mu(X), \frac{M+1}{M-1}A(X)\right)\right)} \,\ge\, \frac{\lambda(\mathcal{B} \cap B(\mu(X),A(X)))}{\left(\frac{M+1}{M-1}\right)^d\lambda(B(\mu(X), A(X)))} \,\ge\,\left(\frac{M-1}{M+1}\right)^d c \,>\, 0\,.$$
This means that conditional on $Y(n) = X$, we may sample from the distribution of $y_{n+1}$ by rejection sampling with success probability bounded away from $0$. That is, we repeatedly sample a point uniformly from $B\left(\mu(X), \frac{M+1}{M-1}A(X)\right)$ until we obtain a sample that lies in $\mathrm{Keep}(X; \mathcal{B})$, and for each trial the probability of success is at least $\left(\frac{M-1}{M+1}\right)^d c$. This property will be used several times. If $\mathcal{B}$ did not have uniform geometry then the success probability for this rejection sampling procedure would not be bounded away from $0$ as $X$ ranges over all sets with sufficiently small~$A(X)$. Fortunately, we may reduce Theorem~\ref{T: original} to the case where $\mathcal{B}$ is bounded and therefore has uniform geometry, as follows.

\begin{lemma} \label{L: reduction to bounded case}
Suppose the conclusion of Theorem~\ref{T: original} holds subject to the extra hypothesis that $\mathcal{B}$ is bounded. Then it also holds without this hypothesis.
\end{lemma}
\begin{proof}
 Let $\mathcal{B}$ be an unbounded convex body and consider the $\mathcal{B}$-valued Jante's law process $Y(n)_{n \in \mathbb{N}_0}$ started at $Y(0) = S$ where $S$ is a set of $M$ distinct points of $\mathcal{B}$. For any $R > 0$, let $\mathcal{B}_R = \mathcal{B} \cap B(0,R)$. Define the random variable $$K := \inf\{ R \in \mathbb{N} \,| \,(\forall n)( Y(n) \subset \mathcal{B}_R)\}\,.$$
 Then $K < \infty$ (and the inf is a min) a.s.,  by Lemma~\ref{L: exponential shrinking} and Proposition~\ref{prop:conv}.
 Let $C = \left\lceil\frac{M+1}{M-1}\sqrt{\frac{(M-1)F(S)}{M^2}} \right\rceil$. Then by Corollary~\ref{cor: F decreases} and \eqref{useful}, we have $\frac{M-1}{M+1}A(Y(n)) \le C$ for all $n \ge 0$. Suppose $R \in \mathbb{N}$ satisfies $\mathbb{P}(K = R) > 0$. Let $\tilde{Y}(\cdot)$ be the $\mathcal{B}_{R+C}$-valued Jante's law process started at $\tilde{Y}(0) = S$, and define $\tilde{K} = \inf\{R \in \mathbb{N} \,|\, (\forall n)( \tilde{Y}(n) \subset \mathcal{B}_{R})\}$. We claim that 
 $\mathbb{P}(K=R) = \mathbb{P}(\tilde{K} = R) $ and that the distribution of $Y(\cdot)$ conditioned on $K = R$ is the same as the distribution of $\tilde{Y}(\cdot)$ conditioned on $\tilde{K} = R$.  Indeed, $Y(\cdot)$ and $\tilde{Y}(\cdot)$ may be coupled by a \emph{maximal coupling}, meaning that for each $n \ge 0$, on the event $Y(n) = \tilde{Y}(n)$, conditional on $Y(n)$ we have $Y(n+1) = \tilde{Y}(n+1)$ with the largest possible probability. When $Y(n) = \tilde{Y}(n)$ and $\mathrm{Keep}(Y(n); \mathcal{B}) = \mathrm{Keep}(\tilde{Y}(n); \mathcal{B}_{R+C})$, this maximal coupling probability is $1$. On the event that $K = R$, this occurs for every $n \ge 0$. For then $\mu(Y(n)) \in \mathcal{B}_R$ for every $n \ge 0$ (since $\mathcal{B}_R$ is convex), and therefore $B(\mu(Y(n)),C) \subseteq B(0,R+C)$ and so
 \begin{eqnarray*}\mathrm{Keep}(Y(n); \mathcal{B})  & = & \mathcal{B} \cap \mathrm{Keep}(Y(n); \mathbb{R}^d) \\ & \subseteq &  \mathcal{B} \cap B\left(\mu(Y(n)), \frac{M+1}{M-1}A(Y(n))\right)\\  & \subseteq & \mathcal{B} \cap B(\mu(Y(n)),C)\\ & \subseteq & \mathcal{B}_{R+C}\,.\end{eqnarray*}
 It follows that for every $n \ge 0$
 $$\mathrm{Keep}(\tilde{Y}(n); \mathcal{B}_{R+C})   =  \mathcal{B}_{R+C} \cap \mathrm{Keep}(\tilde{Y}(n); \mathcal{B})  =  \mathrm{Keep}(\tilde{Y}(n); \mathcal{B}) = \mathrm{Keep}(Y(n); \mathcal{B})$$
as required.

Let $\tilde{\xi}$ be the limit point of $\tilde{Y}(\cdot)$. By the hypothesis that Theorem~\ref{T: original} holds for bounded convex bodies, the distribution of $\tilde{\xi}$ conditioned on $\tilde{K}=R$ is absolutely continuous: it is an absolutely continuous distribution conditioned on an event of positive probability.  Now the distribution of~$\xi$ is a mixture of these conditioned distributions over the set of $R \in \mathbb{N}$ for which $\mathbb{P}(K=R) > 0$. Therefore the distribution of $\xi$ is indeed absolutely continuous (please see also Appendix 1).
\end{proof}

\section{All original points are eventually removed, a.\ s.}\label{S: exodus}
In this section, we will study the $\mathbb{R}^d$-valued Jante's law process $Z(n)_{n \ge 0}$.  In particular, this is an instance of a $\mathcal{B}$-valued Jante's law process where $\mathcal{B}$ has uniform geometry. For some of the lemmas in this section, it takes very little extra work to make the proofs apply to $Y(n)_{n \ge 0}$ subject to the uniform geometry assumption, and we will state them in that generality even though we will only apply them to the process $Z(\cdot)$. So for the whole of this section, we will assume that~$Y(\cdot)$ is a $\mathcal{B}$-valued Jante's law process where $\mathcal{B}$ is a convex body with uniform geometry and  $c = c(\mathcal{B},r_0) >0$.

Let $\tau$ be the first time when all the original points of the configuration are removed, i.e.
\begin{equation}\label{deftau}
\tau=\inf\{n>0:\ Z(n)\cap Z(0)=\emptyset\}\,,    
\end{equation}
with the usual convention that $\tau = \infty$ if there is no $n$ for which $Z(n) \cap Z(0) = \emptyset$.  We know from Proposition~\ref{prop:conv} that for every $\epsilon > 0$ there exists a random $n(\epsilon)$ such that for all $n \ge n(\epsilon)$ we have $Z(n) \subset B(\xi,\epsilon)$, so if for some deterministic choice of $Z(0)$ it happened that $\mathbb{P}(\tau = \infty) > 0$ then we would have $\mathbb{P}(\xi \in Z(0)) > 0$, and the distribution of $\xi$ would not be absolutely continuous. Therefore Theorem~\ref{thm: scale-free limit is continuous} can only hold if it is true that $\tau < \infty$ a.s. It will be useful to prove this fact separately first.  
\begin{prop}\label{prop: exodus}
$\tau<\infty$ a.s.
\end{prop}
The rest of this section is devoted to the proof of Proposition~\ref{prop: exodus}.

In the case $N=3$, i.e.~$M=2$, a.s.~for each $n \ge 0$ the set $Z(n)$ consists of two distinct points, which are equally likely to be removed at the next step. Thus in this simple case, $\tau-1$ is a geometric random variable, in particular a.s.\ finite. 

The case $M=2$ is dealt with, so from now on we will assume that $M \ge 3$. Let $V(d)$ be the volume of the unit ball in $\mathbb{R}^d$. 

Let $\F_n$ denote the $\sigma$-algebra generated by the process $Y(i)_{i \le n}$. Recall from Corollary~\ref{cor: significant decrease of F} that 
\begin{equation}\label{E: increment of log F} 
\mathbb{E}(\log F(Y(n+1)) - \log F(Y(n)) \,|\, \F_n) \le -\frac{4^{-1-d}}{M}\,.
\end{equation}
To complement this we need to bound the downward drift of the smallest inter-point distance~$d(Y(n))$. 
\begin{lemma}\label{lemDhigher}
Suppose that $X$ is a set of $M$ distinct points of $\mathcal{B}$ such that $A(X) \le r_0$. Then
$$
\E(\log d(Y(n)) -\log d(Y(n+1)) \mid Y(n) = X)\le \frac{M}{c d} \left(\frac{d(X)}{A(X)}\right)^d\,.
$$
\end{lemma}
\begin{proof}
Suppose $Y(n) = X$, where $A(X) \le r_0$. The only way the smallest inter-point distance can decrease is for the newly sampled point~$y_{n+1}$ to lie closer than $d(Y(n))$ to one of the existing~$M$ points.  By Lemma~\ref{L: ball inside Keep}, together  with~\eqref{E: uniform geometry 2}, 
$$
\lambda(\mathrm{Keep}(X; \mathcal{B})) \ge \lambda(\mathcal{B} \cap B(\mu(X),A(X)) )\ge  \frac{b(r_0) A(X)^d}{r_0^d} = c V(d) A(X)^d\,.
$$ 
Therefore
\begin{eqnarray*}
\E\left(\left. \log \frac{d(Y(n))}{d(Y(n+1))} \,\right| Y(n) = X \right) & \le &
\frac{M}{\lambda(\mathrm{Keep}(X; \mathcal{B}))} \int_0^{d(X)}  V(d)\, d\, u^{d-1} \log \frac{d(X)}{u} \dd u \\ & \le & \frac{M }{cd} \left(\frac{d(X)}{A(X)}\right)^d.
\end{eqnarray*}
\end{proof}

Next, for any set $X$ of $M$ distinct points in $\mathbb{R}^d$, define
$$
h(X):= \log\left( \frac{\sqrt{F(X)}}{d(X)}\right)=-\log d(X)+\frac{\log F(X)}2\,.
$$
Note that by~\eqref{useful}, we have 
$$
h(X) \ge \frac{1}{2}\log\left(\frac{M(M-1)}{2}\right) \ge \frac12 \log 3> 0\,.
$$
Let $S(d,M)$ be the space of subsets of $\mathbb{R}^d$ of cardinality $M$, with the obvious topology\footnote{ i.e. the topology which it inherits as a quotient of $(\mathbb{R}^d)^M\setminus \Delta$ by the permutation action of the symmetric group $S_M$, where $\Delta \subset (\mathbb{R}^d)^M$ is the subset consisting of sequences in which two or more terms coincide. Note that $S(d,M)$ is not complete with respect to the Hausdorff metric.}. For any $X \in S(d,M)$, denote by~$\hat{X}$ the recentred and rescaled version of $X$, defined by
\begin{equation}\label{hatX}
\hat{X} = \left\{ \frac{x-\mu(X)}{\sqrt{F(X)}}  : x \in X\right\}\,.    
\end{equation}
We have $\mu(\hat{X}) = 0$, $F(\hat{X}) = 1$, and at the same time $h(\hat{X}) = h(X)$. For any $\rho \in [0,\infty)$, define 
\begin{align*}
S(d,M,\rho)&:=\{X \in S(d,M)\,:\, h(X) \le \rho\}\, ,
\\
\hat{S}(d,M,\rho)&:=\{X \in S(d,M)\,:\;h(X) \le \rho, \;\mu(X) = 0,\; F(X) = 1\}\,.
\end{align*}
Then $S(d,M,\rho)$ is invariant under translation and scaling, while $\hat{S}(d,M,\rho)$ is a compact subset of $S(d,M)$.  

Our next goal is to show that 
%for some choice $\rho = \rho_1 < \infty$, 
for the $\mathcal{B}$-valued Jante's law process $Y(\cdot)$, there is a constant $\rho_2 < \infty$ such that a.s. the rescaled and recentred set $\widehat{Y(n)}$ visits the compact set $\hat{S}(d,M,\rho_2)$ infinitely often. 
The proof uses $h$ as a Lyapunov function.

\begin{prop}\label{proppush}
There exist $\rho_1=\rho_1(d,M) < \infty$ and $\gamma_i=\gamma_i(d,M)>0$, $i=1,2$, such that
on the event that $h(Y(n)) \ge \rho_1$ and $A(Y(n)) \le r_0$, we have
\begin{itemize}
\item[(a)] $\E(h(Y(n+1))-h(Y(n))\mid \F_n)\le 0$, and
\item[(b)] $\P(h(Y(n+1))-h(Y(n))<-\g_1 \,|\,\F_n )\ge \g_2$.  
\end{itemize}
\end{prop}
\begin{proof}
For part (a) we use inequality~\eqref{E: increment of log F} and Lemma~\ref{lemDhigher}, along with~\eqref{useful} to compare $A(Y(n))$ and $\sqrt{F(Y(n))}$.  This requires us to take $\rho_1 \ge \frac{1}{d} \left(\log\left(2^{2d+3} M^{d+2}\right) - \log(cd)\right)$.
\begin{figure}
\begin{center}
\includegraphics[scale=0.3]{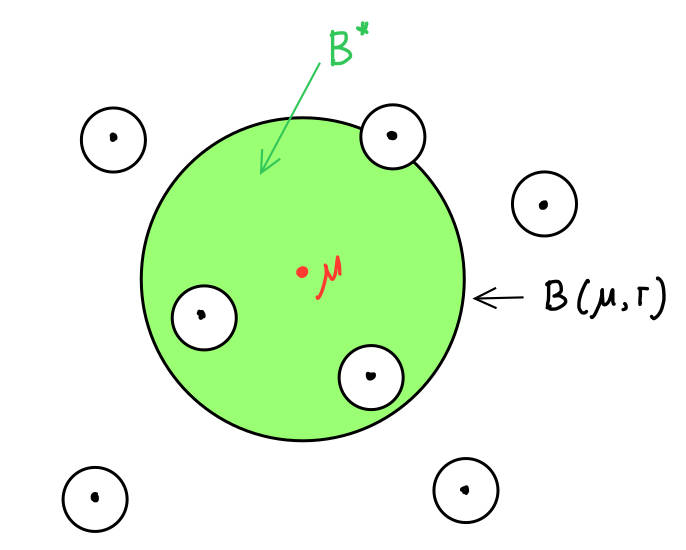}
\end{center}
\caption{A ball around the centre of mass with possible holes}\label{fig-b}
\end{figure}

For part (b),  consider the ball ``with holes" (see Figure~\ref{fig-b})
$$
B^*:= \mathcal{B} \cap \left(B\left(\mu(Y(n)),r \right)
\setminus \bigcup_{w\in Y(n)} 
B\left(w,\left(\frac{c}{2M}\right)^{1/d} r \right)\right)\,,$$
where $$r:=\frac{A(Y(n))}{2}$$ and $c$ is the constant defined in~\eqref{E: c}, which depends on the geometry of $\mathcal{B}$. 

Recall that $Y(n+1)\setminus Y(n)  = \{y_{n+1}\}$. By Lemma~\ref{L: good decrease in F} (with $\alpha = 1/2$), if  $y_{n+1}\in B^*$, then
$$
F(Y(n+1))\le \left(1-\frac{1}{2M}\right)F(Y(n))\,.
$$
On the other hand, using  Lemmas~\ref{L: ball containing Keep} and~\ref{L: ball inside Keep} with Lemma~\ref{L: volume comparison},
$$
\P\left(y_{n+1}\in B\left(\mu(Y(n)),r \right) \right) \ge \frac{\lambda(\mathcal{B} \cap B(\mu(Y(n)),r))}{\lambda\left(\mathcal{B} \cap B\left(\mu(Y(n)),\frac{M+1}{M-1}A(Y(n))\right)\right)}  \ge \frac{(M-1)^{d}}{2^d (M+1)^{d}}\,.
$$  
At the same time, the relative volume of the ``holes" in $B^*$ is bounded. We have $\mu(Y(n)) \in \mathcal{B}$ and $r < r_0$, so we may use the hypothesis of uniform geometry to compare $\lambda(\mathcal{B} \cap B\left(\mu(Y(n)),r \right))$ with $\lambda(B(\mu(Y(n),r)))$. We find
$$
\frac{\lambda(B^*)}{\lambda(B\left(\mu(Y(n)),r \right)\cap \mathcal{B})}
\ge 1- \frac{M}{c} \left(\frac{c}{2M}\right)
= \frac12.
$$
Hence, $\P(y_{n+1}\in B^*)\ge \gamma_2 :=\frac{(M-1)^{d}}{2^{d+1} (M+1)^{d}}>0$. Next, on the event that $y_{n+1}\in B^*$,
$$
\min_{w\in Y(n)} \Vert y_{n+1}-w\Vert\ge \left(\frac{c}{2M}\right)^{1/d} r
\ge \left(\frac{c}{2M}\right)^{1/d}\frac{1}{2M}\sqrt{F(Y(n))}\,.
$$
by~\eqref{useful}. Let $C:= \left(\frac{c}{2M}\right)^{1/d}\frac{1}{2M}$. Then
$$
d(Y(n+1))\ge \min\left\{d(Y(n)),C\sqrt{F(Y(n))}\right\}
$$
and thus
\begin{align*}
h(Y(n+1))-h(Y(n))&=\frac12\log \frac{F(Y(n+1))}{F(Y(n))}-
\log d(Y(n+1))+\log d(Y(n))
\\
&\le  -\frac{1}{4M}-
\min\left\{0,\log \frac{C\sqrt{F(Y(n))}}{d(Y(n))}
\right\}
\\ &
=-\frac{1}{4M}+\max\left\{0,\log \frac{1/C}{h(Y(n))} \right\}
\end{align*}
which equals  $-\g_1:=-\frac1{4M}<0$ provided $h(Y(n))\ge 1/C$.
We see that both parts of the lemma hold when we take
$$ \rho_1 = \max\left(\frac{1}{C}, \frac{1}{d}\log\left(\frac{M^{d+2}4^{d+1}}{cd}\right)\right)
\,.$$
\end{proof}

Let $\rho_2 := \frac{M\,e^{-\rho_1}}{2(M-1)} $ (where $\rho_1$ is the constant from Proposition~\ref{proppush}). Since $e^{-\rho_1}<1$, and  $\frac{M}{2(M-1)}\le 1$ as $M\ge 2$, we get $\rho_2<1$.
Then define
$$
\tilde S(d,M,\rho_2):=\left\{X\in S(d,M) :\ d(X)
\ge \rho_2 \, D(X)\right\}\,.
$$

\begin{prop}\label{propreturns}
The event $\{Y(n)\in\tilde S(d,M,\rho_2)\}$ occurs infinitely often a.s.
\end{prop}
\begin{proof}
First of all, by~\eqref{useful},
$$
X\in S(d,M,\rho_1) \implies X\in \tilde S(d,M,\rho_2)\,.
$$
Indeed, by~\eqref{useful}, we get that $\sqrt{F(X)}\ge \frac{M}{M-1}A(X)\ge\frac{M}{2(M-1)}D(X)$, while $$X\in S(d,M,\rho_1) \Longleftrightarrow h(X)\le \rho_1\Longleftrightarrow  d(X)\ge e^{-\rho_1}\sqrt{F(X)}\,.$$  
For each $m\ge 1$ let $\tau_{m}=\inf\{n\ge m: \ Y(n)\in S(d,M,\rho_1)\}$. By Proposition~\ref{proppush} we have that $\eta_{n;m}=h(Y(n\wedge \tau_m))$ is a non-negative supermartingale for $n\ge m$, for any fixed $m$.  From the supermartingale convergence theorem, it follows that a.s.~$\eta_{n;m}$ converges as $n\to\infty$. Since
%~$\eta_{n;m}$ cannot converge to zero\footnote{in fact, $F(Y(n))\ge \frac{M(M-1)}2 d(X_n)^2$  so that $h(X_n)$ is separated from zero.} and
 on $\{\tau_m>n\}$ we have $\P\left(h(Y(n+1))\le h(Y(n)\right)-\gamma_1)>\gamma_2$, where $\gamma_i>0$, $i=1,2$, are from Proposition~\ref{proppush},  we conclude that a.s.~$\tau_m<\infty$.
As a result, $Y(n)$   returns to $S(d,M,\rho_1)$ and hence also to $\tilde S(d,M,\rho_2)$ infinitely often, a.s. 
\end{proof}

For the remainder of this section, we focus on the process $Z(n)_{n \ge 0}$.

\begin{lemma}\label{lem: points can be removed}
For every set $X$ of $M$ distinct points in~$\mathbb{R}^d$, there exists a finite~$n_0$ and a possible finite trajectory of the process $Z(\cdot)$, say $Z^*(0), \dots, Z^*(n_0)$, such that $Z^*(0) = X$, $Z^*(n_0) \cap X = \emptyset$, along the trajectory there are no ties, i.e. at each step there is a unique legal choice of point~$s_n^*$ to remove, $s_n^* \in Z^*(n-1)$, and the added points $z_1^*, \dots z_n^*$ are distinct and are not in~$X$.
\end{lemma}
\begin{proof} We will identify a trajectory meeting the conditions of the lemma through a combination of probabilistic and constructive reasoning. Note that the process $Z(\cdot)$, started at $Z(0) = X$, a.s.~obeys $F(Z_n) \to 0$, implying that $D(Z_n)\to 0$. Let $\tau_A$ be the least $n$ such that $D(Z_n) < d(X)$. Then $\tau_A < \infty$ a.s., and a.s.~the points $z_1, \dots, z_{\tau_A}$ are distinct and not in $X$. Therefore there exists some $m \in \mathbb{N}$ and a trajectory $Z_0^*, \dots, Z_m^*$, with $Z_i^*\setminus Z_{i-1}^* = \{z_i^*\}$ and $Z_{i-1}^*\setminus Z_{i}^* = \{s_i^*\}$ for $i=1,\dots,m$, such that $D(Z_m^*) < d(X)$, hence $|Z_m^* \cap X| \le 1$, and the points $z_1^*, \dots, z_m^*$ are distinct and not in $X$.  If $Z_m^* \cap X = \emptyset$ then we are done. So from now on assume that $Z_m^* \cap X = \{x_\spad\}$.

The end of the construction falls into two cases, depending on whether $x_\spad$ is a furthest point in $Z_m^*$ from the point 
$$
\mu_0:=\frac{M\mu(Z_m^*) + x_\spad}{M+1}.
$$

First, suppose we are in case 1, where $x_\spad$ is indeed a furthest point in $Z_m^*$ from the point $\mu_0$. If we had $x_\spad = \mu(Z_m^*)$ then we would also have $x_\spad = \mu_0$ and hence all points of $Z_m^*$ would have to coincide, which was ruled out in the construction of $Z_m^*$. Therefore $x_\spad \neq \mu(Z_m^*)$. Let $z_{m+1}^* = \epsilon \mu(Z_m^*) + (1-\epsilon) x_\spad$, for some very small $\epsilon \in (0,1)$, chosen so that~$z_{m+1}^*$ is not in $\bigcup_{i=0}^m Z_i^*$. (This can fail only for finitely many values of $\epsilon$.) Then $x_\spad$ is the unique furthest point of $Z_m^* \cup \{z_{m+1}^*\}$ from its centre of mass, which is 
$$
\mu_\epsilon:=\frac{(M+\epsilon)\mu(Z_m^*) + (1-\epsilon)x_\spad}{M+1}\,.
$$ 
To see this, note first that $$\Vert z_{m+1}^* - \mu_\epsilon\Vert = \frac{M}{M+1}(1-\epsilon)\Vert x_\spad - \mu(Z_m^*) \Vert > 0\,,$$ while $$\Vert x_\spad - \mu_\epsilon\Vert = \frac{M+\epsilon}{M+1}\Vert x_\spad - \mu(Z_m^*) \Vert\,.$$ Thus $\Vert x_\spad - \mu_\epsilon\Vert > \Vert z_{m+1}^* - \mu_\epsilon\Vert > 0$. Secondly, note that the ball $B(\mu_\epsilon, \Vert x_\spad - \mu_\epsilon\Vert)$  contains the ball $B(\mu_0, \Vert x_\spad - \mu_0\Vert)$, and their boundary spheres meet only at the point $x_\spad$. Hence $\Vert x_\spad - \mu_\epsilon\Vert > \Vert z - \mu_\epsilon\Vert$ for every $z \in Z_m^* \setminus \{x_\spad\}$. Hence the point $x_\spad$ is the unique point which may be removed, and we obtain $Z_{m+1}^* = (Z_{m} \cup \{z_{m+1}^*\}) \setminus \{x_\spad\}$. The trajectory $Z_0^*, \dots, Z_{m+1}^*$ is a permissible trajectory with no tiebreaks and no repeated points, and $Z_{m+1}^* \cap X = \emptyset$, as required. 

Case 2 is the case where $x_\spad$ is not a furthest point in $Z_m^*$ from the point $\mu_0$. For the rest of the proof, we assume that we are in this case. We claim that if $v \in \mathbb{R}^d\setminus\{0\}$ is sufficiently close to $0 \in \mathbb{R}^d$ then from the trajectory $Z_0^*, \dots, Z_m^*$ we can obtain another permissible trajectory in which the arriving points are all perturbed by the vector~$v$. The perturbed trajectory is $Z_0^\dagger, \dots, Z_m^\dagger$, where $Z_0^\dagger = X$ and for $i=1, \dots, m$ we have $Z_i^\dagger = ( Z_{i-1}^\dagger \cup \{z_i^\dagger\})\setminus\{s_i^\dagger\}$, where $z_i^\dagger = z_i^* + v$ and $s_i^\dagger$ is a point of $Z_{i-1}^\dagger \cup \{z_i^\dagger\}$ furthest from $\mu(Z_{i-1}^\dagger \cup \{z_i^\dagger\})$.  To say that the trajectory $Z_0^\dagger, \dots, Z_m^\dagger$ is permissible means that $z_i^\dagger \in \mathrm{Keep}(Z_{i-1}^\dagger; \mathbb{R}^d)$ for each $i = 1, \dots, m$, so that $z_i^\dagger \neq s_i^\dagger$.  In fact, if~$v$ is sufficiently close to~$0$ then the corresponding points are removed at each step, by which we mean that $Z_i^\dagger = Z_{i-1}^\dagger \cup \{z_i^\dagger\} \setminus \{s_i^\dagger\}$, where either $s_i^\dagger = s_i^*$ (if $s_i \in X$), or $s_i^\dagger = s_i^* + v$ (if $s_i \notin X$).  To see this, note that each point $z_i^\dagger$ or $s_i^\dagger$ depends continuously on $v$, as does $\mu(Z_i^\dagger)$. The (assumed) statement that $s_i^*$ is the unique point of $Z_{i-1}^* \cup \{z_i^*\}$ furthest from $\mu(Z_{i-1}^* \cup \{z_i^*\})$ is equivalent to a finite collection of strict inequalities:  for each $z \in (Z_i^* \cup \{z_i^*\})\setminus\{s_i^*\}$, $$ \Vert s_i^* - \mu(Z_i^* \cup \{z_i^*\}) \Vert > \Vert z - \mu(Z_i^* \cup \{z_i^*\})\Vert\,. $$ It follows that when $v \in \mathbb{R}^d \setminus\{0\}$ is sufficiently close to $0$ then for each $z \in (Z_{i-1}^\dagger \cup \{z_i^\dagger\})\setminus\{s_i^\dagger\}$, $$ \Vert s_i^\dagger - \mu(Z_i^\dagger \cup \{z_i^\dagger\}) \Vert > \Vert z - \mu(Z_i^\dagger \cup \{z_i^\dagger\})\Vert\,,$$ i.e.\ $s_i^\dagger$ is the unique point of $Z_{i-1}^\dagger \cup \{z_i^\dagger\}$ furthest from $\mu(Z_{i-1}^\dagger \cup \{z_i^\dagger\})$ for each $i=1, \dots, m$. In particular, the trajectory $Z_0^\dagger, \dots, Z_m^\dagger$ involves no tiebreaks when $v$ is sufficiently close to $0$. 
By similar reasoning, since the trajectory $Z_0^*, \dots, Z_m^*$ involves no repeated points, the same is true of the trajectory $Z_0^\dagger, \dots, Z_m^\dagger$ when $v$ is sufficiently close to $0$.

Having chosen a vector $v$ suitably close to $0$, we now extend our two trajectories one further step, so that $Z_{m+1}^\dagger = Z_m^* + v$, i.e.\ $Z_m^\dagger$ is the translate of $Z_m^*$ by the vector $v$. To do this, we choose $z_{m+1}^* = x_\spad -v$ and $z_{m+1}^\dagger = x_\spad + v$. 
Then 
\begin{equation}\label{E: happy translation} 
Z_{m}^\dagger \cup \{z_{m+1}^\dagger\} = \left( Z_{m}^* \cup \{z_{m+1}^*\}\right) + v\,.
\end{equation} 
Notice that in the translation which relates the two configurations, the roles of $x_\spad$ and the new point ($z_{m+1}^*$ and $z_{m+1}^\dagger$ respectively) are swapped. We may choose $v$ so that there is no tiebreak in either configuration when selecting the point furthest from the centre of mass. This is because the set of vectors $v$ which would cause such a tie is contained in a finite union of hypersurfaces of dimension less than $d$, so we may choose $v$ as close as we like to $0$ but not in any of these hypersurfaces. From \eqref{E: happy translation}, we see that there are points $s_{m+1}^* \in Z_{m}^* \cup \{z_{m+1}^*\}$ and $s_{m+1}^\dagger \in Z_{m}^\dagger \cup \{z_{m+1}^\dagger\}$ which are furthest from the centres of mass $\mu(Z_{m}^* \cup \{z_{m+1}^*\})$ and $\mu(Z_{m}^\dagger \cup \{z_{m+1}^\dagger\})$ respectively, and which are related by $s_{m+1}^\dagger  = s_{m+1}^* + v$. Hence $Z_{m+1}^\dagger = Z_{m+1}^* + v$. Before proceeding, we must check that the trajectories $Z_0^*, \dots, Z_{m+1}^*$ and $Z_0^\dagger, \dots, Z_{m+1}^\dagger$ are permissible. This is where we use the assumption that we are in case 2, together with the continuity of all defined points as functions of $v$.  Let $u^*$ be a point of $Z_m^*$ which is strictly further from $\mu_0$ than $x_\spad$ is from $\mu_0$, and let $u^\dagger = u^* + v$. As $v \to 0$, we have $\mu(Z_{m}^* \cup \{z_{m+1}^*\}) \to \mu_0$ and $\mu(Z_{m}^\dagger \cup \{z_{m+1}^\dagger\}) \to \mu_0$, while $z_{m+1}^* \to x_\spad$, $z_{m+1}^\dagger \to x_\spad$ and $u^\dagger \to u^*$. Therefore for $v$ close enough to $0$ we have $$\Vert u^* - \mu(Z_{m}^* \cup \{z_{m+1}^*\})\Vert > \Vert z_{m+1}^* - \mu(Z_{m}^* \cup \{z_{m+1}^*\})\Vert $$ and $$\Vert u^\dagger - \mu(Z_{m}^\dagger \cup \{z_{m+1}^\dagger\})\Vert > \Vert z_{m+1}^\dagger - \mu(Z_{m}^\dagger \cup \{z_{m+1}^\dagger\})\Vert\,, $$ showing that both trajectories are indeed permissible. 

Finally, by the same argument as in the first paragraph of this proof, there exists a possible trajectory $Z_{m+1}^*, \dots, Z_{m+m'}^*$ such that $|Z_{m+1}^* \cap Z_{m+m'}^*| \le 1$, and having no tiebreaks, no repeated points, and no arriving points which were already in $\bigcup_{i=0}^m Z_i^*$.  There are now two possibilities: either $x_\spad \notin Z_{m+m'}^*$, in which case we are done, or $x_\spad \notin Z_{m+m'}^*$. In the latter case, consider the trajectory $Z_{m+1}^\dagger, \dots, Z_{m+m'}^\dagger$ defined by $Z_{i}^\dagger = Z_i^* + v$ for $i = m+1, \dots, m+m'$. It is also a permissible trajectory, being a translation of a permissible trajectory. Moreover, $Z_{m+1}^\dagger \cap Z_{m+m'}^\dagger = \{x_\spad + v\}$, and in particular $x_\spad \notin Z_{m+m'}^\dagger$, so that $Z_{m+m'}^\dagger \cap X = \emptyset$. Therefore $Z_0^\dagger, \dots, Z_{m+m'}^\dagger$ is a trajectory satisfying the conclusion of the lemma, and we are done.

\end{proof}

Recall the following definition.
\begin{defn}
Let $X$ and $X'$ be two subsets of $\R^d$. The {\em Hausdorff distance} between them is defined as
\begin{align*}
d_H(X,X')&=\max\left\{
\sup_{x\in X} \left[\inf_{a'\in X'} \Vert x-a'\Vert\right],
\sup_{x'\in X'} \left[\inf_{a\in X} \Vert x'-a\Vert\right]
\right\}
\\ & =
\inf\left\{\epsilon>0:\ \displaystyle X\subseteq 
\bigcup_{a'\in X'} \mathrm{B}(a',\epsilon)
\text{ and } X'\subseteq 
\bigcup_{a\in X} \mathrm{B}(a,\epsilon)
\right\}.
\end{align*}
\end{defn}
Note that when $X$ and $X'$ are both finite, then $d_H(X,X')=0$ if and only if  $X$ and $X'$ have the same closure.

\begin{lemma}\label{lem: neighbourhood for point removal}
For every set $X \in S(d,M)$, there is a finite $n_0 = n_0(X)$, an~$\epsilon_0 = \epsilon_0(X) \in (0, d(X)/2)$, and a  $\delta_0 = \delta_0(X) > 0$ such that if $X' \in S(d,M)$ has Hausdorff distance  less than~$\epsilon_0$ from $X$ then 
 $$
 \mathbb{P}(Z(n_0) \cap Z(0) = \emptyset \mid Z(0) = X')  > \delta_0\,.
 $$ 
\end{lemma}
\begin{proof}
Let $z_{1-M}^*, \dots, z_0^*$ be any $M$ distinct points in $\mathbb{R}^d$, and let $X = \{z_{1-M}^*, \dots, z_0^*\}$. Take a trajectory as provided by Lemma~\ref{lem: points can be removed}.   Since at each step $1 \le n \le n_0$ there is a unique choice of point $s_n^* \in Z^*(n-1) \cup \{z^*_{n}\}$ such that $s_n^*$ maximizes the distance from $\mu(Z^*(n-1) \cup \{z_n^*\})$, and there are only finitely many steps, there is an $\epsilon > 0$ such that for all $1 \le n \le n_0$ we have for any $z \in Z^*(n)$ that $$
\|s_n^* - \mu(Z^*(n-1) \cup \{z_n^*\})\| \ge \|z - \mu(Z^*(n-1) \cup \{z_n^*\})\| + 4 \epsilon\,.
$$ 
By making $\epsilon$ smaller if necessary, we may also assume that 
$$
\min(\|z-x\|\,:\,z \in Z^*(n_0),\, x \in X) > 2\epsilon.
$$
Consider any perturbed trajectory $Z'(0), \dots, Z'(n_0)$ where $Z_0' = \{z_{1-M}', \dots, z_0'\}$ and for $1 \le n \le n_0$ we have $Z'(n) = (\{z_n'\} \cup Z'(n-1)) \setminus \{s_n'\}$, and for each $i = 1-M, \dots, n_0$ we have $\|z_i^* - z_i'\| < \epsilon$.
Then one may prove by induction over~$n$ that in the perturbed trajectory, the thrown out points~$s_n'$ correspond to the original points~$s_n$ in the sense that for each $1 \le n \le n_0$, the index~$j$ such that $s_n^* = z_j^*$ also satisfies $s_n' = z_j'$. This is because the perturbation moves each point by a distance less than~$\epsilon$ and also moves each centre of mass by a distance at most~$\epsilon$. It therefore alters corresponding distances by less than~$2\epsilon$, and differences of distances by less than~$4\epsilon$. 

In particular, the points of $Z'(n_0)$ are in one-to-one correspondence with points of $Z^*(n_0)$ so that corresponding pairs of points are at a distance less than $\epsilon$. It follows that $Z'(n_0) \cap X' = \emptyset$.

Finally, consider the random process $Z(\cdot)$ started at $Z(0) = X'$. The probability that each new point $z_i$, $i = 1, \dots, n_0$, lies in the ball $B(z_i^*,\epsilon)$ is bounded away from $0$, uniformly over all choices of $X'$ that are $\epsilon$-perturbations of $X$. This is because the volumes of the sets $\mathrm{Keep}(Z(n);\mathbb{R}^d)$ are bounded above \emph{a priori} in terms of $X$, using the fact that $$\frac{M^2}{M-1}A(Z_n)^2 \le F(Z_n) \le F(X') \le \frac{M(M-1)}{2}D(X')^2 \le \frac{M(M-1)}{2}(D(X) + 2\epsilon)^2$$ together with Lemma~\ref{L: ball containing Keep}. When this occurs, we have $Z(n_0) \cap Z(0) = \emptyset$.
\end{proof}

\begin{lemma}\label{lem: point removal on the good set}
There exist $n_1 = n_1(d,M) < \infty$ and $\delta_1 = \delta_1(d,M) > 0$ such that for every $X \in \tilde{S}(d,M,\rho_2)$,
$$\mathbb{P}(Z(n+n_1) \cap Z(n) = \emptyset \mid  Z(n) = X) \ge \delta_1\,.$$
\end{lemma}
\begin{proof}
First, observe that since the process $Z(\cdot)$ is time-homogeneous and invariant under translation and scaling, we have
$$
\mathbb{P}(Z(n+n_1) \cap Z(n) = \emptyset \mid  Z(n) = X) = \mathbb{P}(Z(n_1) \cap Z(0) = \emptyset \mid  Z(0) = \hat{X})\,,
$$
where $\hat{X}$ is defined by~\eqref{hatX}.
So it suffices to prove that there exist $n_1 < \infty$ and $\delta_1 > 0$ such that for all $X \in \hat{S}(d,M,\rho_2)$ we have
$$
\mathbb{P}(Z(n_1) \cap Z(0) = \emptyset \mid Z(0) = X) \ge \delta_1\,.
$$
This follows from Lemma~\ref{lem: neighbourhood for point removal}, since for each point $X$ in the compact set $\hat{S}(d,M,\rho_2)$, that statement specifies a neighbourhood $U_X$ of $X$ and constants $n_0(X) < \infty$, $\delta_0(X) > 0$ such that for all $X' \in U_X$, 
$$
\mathbb{P}(Z(n_0) \cap Z(0) = \emptyset \mid\, Z(0) = X') \ge \delta_0(X)\,.
$$ 
Since $\hat{S}(d,M,\rho_2)$ is compact, there exists a covering of $\hat{S}(d,M,\rho_2)$ by some finite sequence of such neighbourhoods $U_{X_1}, \dots, U_{X_k}$, and we may take $$n_1 = \max(n_0(X_1), \dots, n_0(X_k)) < \infty$$ and $$ \delta_1 = \min(\delta_0(X_1), \dots, \delta_0(X_k)) > 0\,.$$ 
\end{proof}

We can now prove Proposition~\ref{prop: exodus}, which claims that $\tau < \infty$ a.\ s.   
\begin{proof}[Proof of Proposition~\ref{prop: exodus}]
 Starting from any configuration $Z(0) = X$ of $M$ distinct points in $\mathbb{R}^d$, the process $Z(\cdot)$ a.s.~eventually enters $\tilde{S}(d,M,\rho_2)$, by Proposition~\ref{propreturns}, say at stopping time $\tau_1$. Let $n_1 = n_1(d,M)$ be the constant provided by Lemma~\ref{lem: point removal on the good set}. Inductively define $\tau_k$ for $k=2,3,\dots$, by
 $$ \tau_k = \min\{ t\,:\, t > \tau_{k-1}+n_1 , \, Z(t) \in \tilde{S}(d,M,\rho_2)\}\,.$$
 For each $k \ge 1$, conditional on $\tau_k$ and $Z(\tau_k)$, in time steps $\tau_k +1, \dots, \tau_k + n_1$ all the points of $Z(\tau_k)$  are removed, with probability at least $\delta_1$, regardless of the configuration $Z(\tau_k)$. Therefore with probability $1$ there exists a finite $k$ for which this event occurs. We then have $\tau < \tau_k + n_1 < \infty$, as required. \end{proof}

\section{Proof of Theorem 1}\label{S: continuous limit for scale-free process}

Without loss of generality, we assume that the initial state $Z(0)$ is deterministic. This is harmless since if for each deterministic choice of $Z(0)$ the limit $z_\infty$ exists a.s.~and has an absolutely continuous distribution, then if instead $Z(0)$ is random, $z_\infty$ still exists a.s.~and its distribution is a mixture of absolutely continuous distributions, which is necessarily absolutely continuous.

List the elements of $Z(0)$ in an arbitrary order as $z_{-(N-2)}, z_{-(N-3)}, \dots z_{0}$. Recalling that for $n \ge 1$, $z_n$ is the unique point of $Z(n) \setminus Z(n-1)$, we observe that for each $m \ge 0$ we have 
$$ 
\{z_i: i = -(N-2),...,m\}\;=\; \bigcup_{n=0}^m Z(n) \,.
$$
We discarded a null set to ensure that the points of the sequence $(z_i)_{i=-N_2}^\infty$ are pairwise distinct. Thus for each $n \ge 1$, the unique point $s_n$ of $Z(n-1) \setminus Z(n)$ is equal to $z_{\a(n)}$ for a unique index $\a(n)$ in the range $-(N-2) \le \a(n) \le n-1$. Moreover, the indices $\a(1), \a(2), \a(3), \dots$ are pairwise distinct random variables since each point is removed at most once.

\begin{lemma}
Fix any $K\ge N-1$ and any deterministic sequence of distinct indices $\mathbf{i} = (i(1),\dots,i(K))$ such that for each $n = 1, \dots, K$ we have $-(N-2) \le i(n) \le n-1$ and $K$ is the least positive integer such that $\{-(N-2), \dots, 0\} \subseteq \{i(1), \dots, i(K)\}$. Let ${\bar\a}=(\a(1),\dots,\a(K))$. Then the event that $\tau = K$ and $\bar{\a} =\mathbf{i}$ is equivalent to the event that a certain finite collection of linear or quadratic inequalities are satisfied, which involve the coordinates of the random points $z_1, \dots, z_K$, as well as the coordinates of the deterministic points $z_{-(N-2)}, \dots, z_0$. 
\end{lemma}
\begin{figure}
\begin{center}
\includegraphics[scale=0.4]{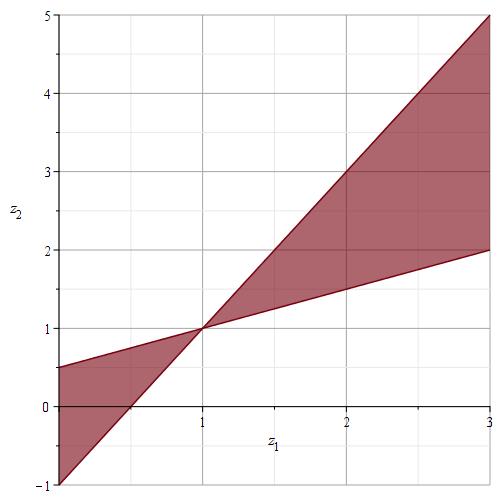}
\end{center}
\caption{Area where $(z_1,z_2)$ can be sampled in the example of Remark~\ref{rem1}}\label{fig2}
\end{figure}

\begin{rema}\label{rem1}
For example, take $d=1$, $N=3$, and w.l.o.g.\ $z_{-1}=-1$, $z_0=1$. Then on $\tau=2$, $\a(1)=-1$, $\a(2)=(0)$, we have 
$z_1\in (0,1)$ and $z_2\in\left(2z_1-1,\frac{z_1+1}{2}\right)$
or
$z_1\in (1,3)$ and $z_2\in\left(\frac{z_1+1}{2},2z_1-1\right)$. See Figure~\ref{fig2}.
\end{rema}

These inequalities define a bounded semi-algebraic set\footnote{i.e., a set defined by a number of polynomial inequalities and equalities; in our case, a.s.~these will be just inequalities.} $\mathcal{P}_{K,\mathbf{i}} \subset \left(\mathbb{R}^d\right)^K$.
(When $d=1$, $\mathcal{P}$ is a polytope, since the defining inequalities may be reduced to a collection of affine linear inequalities.)  For some sequences $\mathbf{i}$, $\mathcal{P}_{K,\mathbf{i}}$ may be empty or have empty interior. 
Note that the boundary of $\mathcal{P}_{K,\mathbf{i}}$  is contained in a finite union of algebraic hypersurfaces where at least one of the defining inequalities becomes an equality. These hypersurfaces have trivial $(Kd)$-dimensional Lebesgue measure. 
Thus $\mathcal{P}_{K,\mathbf{i}}$ has positive $(Kd)$-dimensional Lebesgue measure if and only if it has a non-empty interior.

Let $E_{K,\mathbf{i}}$ be the event that $\tau = K$ and $\bar\a=\mathbf{i}$.  If $\mathbb{P}\left(E_{K,\mathbf{i}}\right) > 0$ then conditional on $E_{K,\mathbf{i}}$ the sequence $\bar\a$ is distributed according to $U\left(\mathcal{P}_{K,\mathbf{i}}\right)$; 
in particular it lies a.s.~in the interior~$\mathcal{P}_{K,\mathbf{i}}^\circ$.

  Define a $d$-dimensional linear subspace $L_K$ of $\left(\mathbb{R}^d\right)^K$ by 
$$L_K  := \{ \left((v, \dots, v) \in \mathbb{R}^d\right)^K : v \in \mathbb{R}^d\}\,.$$ 
Consider the orthogonal projection $\pi_K$ of $\left(\mathbb{R}^d\right)^K$ onto $\left(L_K\right)^\perp$. We have
$$ \lambda_{Kd}(\mathcal{P}_{K,\mathbf{i}}) = \int_{\pi_K(\mathcal{P}_{K,\mathbf{i}})} \lambda_{d}(\mathcal{P}_{K,\mathbf{i}} \cap \pi_K^{-1}(x))\,d\lambda_{(K-1)d}(x)\,. $$
We conclude that if $\lambda_{Kd}\left(\mathcal{P}_{K,\mathbf{i}}\right) > 0$ then for $\lambda_{Kd}$-a.e. $\mathbf{p} = (p_1, \dots, p_K) \in \mathcal{P}_{K,\mathbf{i}}$, the set 
$$ 
T(\mathbf{p},K,\mathbf{i}): = \{v \in \mathbb{R}^d\,:\, (p_1 + v, \dots, p_K + v) \in \mathcal{P}_{K,\mathbf{i}} \} 
$$
has $\lambda_d(T(\mathbf{p},K,\mathbf{i})) > 0$. 
Therefore we can sample a $U(\mathcal{P}_{K,\mathbf{i}})$ random variable $\mathbf{q}$ by sampling first a $U(\mathcal{P}_{K,\mathbf{i}})$ random variable $\mathbf{p}$, then (conditional on $\mathbf{p}$) sampling a $U(T(\mathbf{p},K,\mathbf{i}))$ random variable~$\delta$ and setting $\mathbf{q} := \mathbf{p} + (\underbrace{\delta,\delta, \dots, \delta}_{K\text{ times}})$.

Recall the definition of $\tau$ from~\eqref{deftau} and
let $\mathcal{F}_\tau$ be the $\sigma$-algebra generated by the stopped process $Z(\cdot \wedge \tau)$. Let $V$ be a random variable whose distribution conditional on $\mathcal{F}_\tau$ is $$U(T((z_1, \dots, z_\tau),\tau, (\a(1), \dots, \a(\tau))))\,,$$ such that~$V$ is conditionally independent of $(z_{\tau+1}, z_{\tau+2}, \dots)$ given $\mathcal{F}_\tau$. That is to say,  on the event~$E_{K,\mathbf{i}}$, given $z_{-(N-2)}, \dots, z_K$, we let $V \sim U(T((z_1, \dots, z_K), K,\mathbf{i}))$, independent of $(z_{K+1}, z_{K+2}, \dots)$.  Proposition~\ref{prop: exodus}, and the fact that ties a.s.~do not occur, imply that a.s.~$E_{K,\mathbf{i}}$ holds for exactly one $(K,\mathbf{i})$, so $V$ is well-defined. 
Since $V$ is the uniform distribution on a subset of $\mathbb{R}^d$ of positive and finite $d$-dimensional Lebesgue measure, the conditional distribution of $V$ on $\mathcal{F}_\tau$ is a.s.~absolutely continuous. This is a crucial property which we will use later.
%The crucial property of $V$ is that its conditional distribution of $V$ on $\mathcal{F}_\tau$ is a.s.~continuous, being the uniform distribution on a subset of $\mathbb{R}^d$ of positive and finite $d$-dimensional Lebesgue measure. 

We will now show that we can ``perturb'' the original path of new points $z_1, \dots, z_\tau$ by adding the same absolutely continuous random vector to each of them, without changing the distribution of the path. (Of course, this requires that the vector that we add is not independent of $(z_1, \dots, z_\tau)$.) This implies that the distribution of the set $Z(\tau)$ is a mixture of distributions, each of which is the distribution of the translation by some absolutely continuous random vector of some set of $M$ distinct points. This is the key argument in showing that the limit point $z_\infty$ has an absolutely continuous distribution.

Define a new random sequence $\left(z_n'\right)_{n=-(N-2)}^\infty$ as follows. 
$$
z_n' = \begin{cases} z_n, & \text{if $n \le 0$,}\\ z_n + V, & \text{if $n \ge 1$,} \end{cases}
$$
Conditional on $z_{-(N-2)}, \dots, z_0$ and on $E_{K,\mathbf{i}}$, the sequences $\left(z_n\right)_{n=1}^\tau$ and $\left(z_n'\right)_{n=1}^\tau$ are identically distributed: both are $U(\mathcal{P}_{K,\mathbf{i}})$ random variables.  

Let $Z'(0) = Z(0)$, then for $n \ge 1$ let 
$$
Z'(n) =  \{z_n'\} \cup \left(Z'(n-1) \setminus \{z_{\a(n)}\}\right)\,.
$$
For all $n \ge \tau$ we have 
\begin{equation}\label{E: translation} Z'(n) = \{z + V \,:\, z \in Z(n)\}
\end{equation}
Moreover $Z'(\cdot)$ is a Markov chain with the same distribution as $Z(\cdot)$: by construction they have the same distribution up to the random time $\tau$, and from time $\tau$ onwards, the transitions have the correct distribution by equation~\eqref{E: translation} together with the translation-invariance of the transition law of $Z(\cdot)$.

Therefore the sequence $\left(z_n'\right)$ has the same distribution as the original sequence $(z_n)$.  Since~$\left(z_n\right)$ a.s.~converges to a limit $z_\infty$, it also holds that $\left(z_n'\right)$ a.s.~converges; let $z_\infty'$ be its limit. The distribution of $z_\infty'$ coincides with the distribution of $z_\infty$. 

On the other hand, the distribution of $z_\infty'$ is also the mixture of its conditional distributions on $\mathcal{F}_\tau$. But we have 
$$
z_\infty' = z_\infty + V
$$
and $z_\infty$ and $V$ are conditionally independent given $\mathcal{F}_\tau$. Thus the distribution of $z_\infty'$, conditionally on $\F_\tau$, is a
distribution of a sum of two independent random variables, one of which (namely~$V$) has an absolutely continuous distribution. It is well-known that such a sum is also an absolutely continuous random variable (see e.g.~Theorem 2.1.16 in~\cite{DUR}). Therefore the {\em unconditional} distribution of $z_\infty'$ is also absolutely continuous, being a  mixture of absolutely continuous distributions. This completes the proof of Theorem~\ref{thm: scale-free limit is continuous}.

\section{\texorpdfstring{Coupling $Y(\cdot)$ and $Z(\cdot)$}{Lg}}
\label{S: coupling succeeds}
In this section, we complete the proof of Theorem~\ref{T: original}. By Lemma~\ref{L: reduction to bounded case} it suffices to prove the absolute continuity of the limit point $\xi$ for a $\mathcal{B}$-valued Jante's law process where $\mathcal{B}$ is bounded. So we shall assume throughout this section that $\mathcal{B}$ is bounded, and therefore has uniform geometry.

\begin{figure}[ht]
\begin{center}
\includegraphics[scale=0.4]{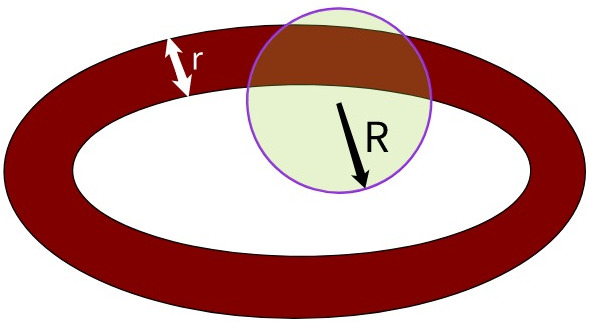}
\end{center}
\caption{Lemma~\ref{L: easy isoperimetric} bounds the quantity of chocolate in the ball of radius $R$.}
\label{F: choco egg}
\end{figure}

We will need an isoperimetric inequality for inner shells of convex bodies, for which we have been unable to find a reference. It concerns the following problem. Suppose you have a (possibly) hollow chocolate egg whose outer boundary is the boundary of a convex body. If all the chocolate is within distance $r$ of the outer boundary of the egg, what is the maximum quantity of chocolate that can possibly be contained within a ball of radius $R$? See Figure~\ref{F: choco egg}.

\begin{lemma}\label{L: easy isoperimetric} Let $R>r >0$ and let $\mathcal{B}$ be a convex body in $\mathbb{R}^d$. For every $y \in \mathbb{R}^d$, we have
 \begin{equation}\label{E: easy isoperimetric} 
 \lambda\left(B(y,R) \cap \mathcal{B} \cap N_r(\mathcal{B}^c)\right)  \le 2rd\sqrt{d}\,V(d-1)R^{d-1}\,,
 \end{equation} where $V(d-1)$ is the volume of the unit ball in $\mathbb{R}^{d-1}$, with the convention that $V(0) = 1$.
\end{lemma}
\begin{proof}
For $d=1$ the inequality holds with $V(0) = 1$, since the set $N_r\left(\mathcal{B}^c\right)$ is a union of at most two intervals of length at most $r$ and therefore has length at most $2r$. Now assume that $d \ge 2$. Let $A = \mathcal{B} \cap B(y,R) \cap N_r(\mathcal{B}^c)$. Consider any point $x \in A$. Then $\Vert x-z\Vert \le r$ for some  $z \in \partial\mathcal{B}$. Let $H$ be a supporting hyperplane through $z$, meaning that $H \cap \mathcal{B}^\circ = \emptyset$. Then (as for any hyperplane through $z$) $H$ cuts at least one of the $d$ axis-parallel lines through $x$ at a point whose distance from $x$ is at most $r\sqrt{d}$.  It follows that there is a point of $\partial\mathcal{B}$ that lies between $x$ and $H$ on one of the axis-parallel lines through $x$. The situation is illustrated in Fig.~\ref{F: supporting hyperplane} below. 

\begin{figure}[ht]
\begin{center}
\includegraphics[scale=0.15]{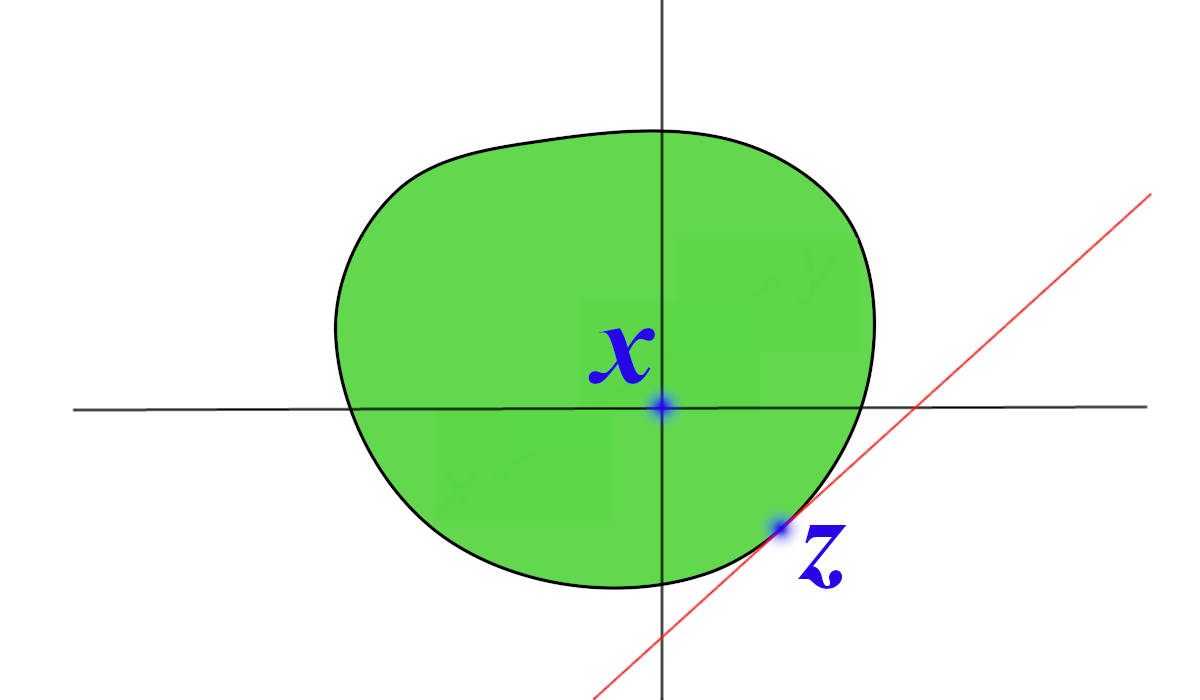}
\end{center}
\caption{The supporting hyperplane at $z$ to the green convex body $\mathcal{B}$ cuts at least one axis-parallel line through $x$ within distance $\sqrt{d}\Vert z-x \Vert$ of $x$.}
\label{F: supporting hyperplane}
\end{figure}

Thus $A$ is covered by the $d$ sets $A_1, ..., A_d$, where 
$$ 
A_i = \left\{ x \in A \,:\, \left(\exists t \in \left[-r\sqrt{d},r\sqrt{d}\right]\right)\left(x + t e_i \in \partial \mathcal{B}\right)  \right\}\,,
$$
where $e_1, \dots, e_d$ are the standard basis vectors of $\mathbb{R}^d$. Since $A_i \subset B(y,R)$, the orthogonal projection of $A_i$ onto a hyperplane orthogonal to $e_i$ has $d-1$-dimensional Lebesgue measure at most $V(d-1)R^{d-1}$. Each line parallel to $e_i$ that meets $B(y,R)$ meets $A_i$ in at most two intervals, whose total length is at most $2r\sqrt{d}$. Hence $\lambda(A_i) \le 2r\sqrt{d} V(d-1)R^{d-1}$ for each $i = 1, \dots, d$. Summing over $i$ we get $$\lambda(A) \le \lambda\left(\bigcup_{i=1}^d A_i\right) \le \sum_{i=1}^d \lambda(A_i) \le 2rd\sqrt{d}V(d-1)R^{d-1}\,.$$
\end{proof}

Lemma~\ref{L: easy isoperimetric} will suffice for our purposes, but we remark that the sharp version of the inequality is as follows. (A proof is given in the appendix.) 
\begin{lemma}\label{L: isoperimetric}
 Let $R>r >0$ and let $\mathcal{B}$ be a convex body in $\mathbb{R}^d$. For every $x \in \mathbb{R}^d$, we have
 \begin{equation}\label{E: isoperimetric} 
 \lambda\left(B(x,R) \cap \mathcal{B} \cap N_r(\mathcal{B}^c)\right)  \le V(d)(R^d - (R-r)^d)\,,
 \end{equation} where $V(d)$ is the volume of the unit ball in $\mathbb{R}^d$.
\end{lemma}

Let $\mathcal{B}\subset \mathbb{R}^d$ be a bounded convex body. For any set $X$ of $M$ distinct points of $\mathcal{B}$ define
$$ \mathrm{d}_{\mathcal{B}}(X) = \mathrm{dist}(X, \mathcal{B}^c)\,.
$$
and for any set $X$ of $M$ distinct points of $\mathcal{B}$ that is not entirely contained in $\partial{\mathcal{B}}$, define
$$\mathrm{d}^{\circ}_{\mathcal{B}}(X) = \mathrm{dist}(X \cap \mathcal{B}^\circ, \mathcal{B}^c)\,,
$$
For $Y(n)_{n \ge 0}$ a $\mathcal{B}$-valued Jante's law process we define
$$
g_n := \begin{cases}\frac12 \log (F(Y(n))) - \log\left(\mathrm{d}^{\circ}_{\mathcal{B}}(Y(n)\right),
&\text{if }Y(n) \not\subset \partial\mathcal{B};
\\
+\infty,&\text{if }Y(n) \subset \partial \mathcal{B}\,.
\end{cases}
$$
\begin{lemma}\label{L: second supermartingale}
There exists $\delta > 0$ (depending on $M$, $d$, and $\mathcal{B}$) such that for any set $X_0$ of $M$ distinct points of $\mathcal{B}$ and any $n_0$, 
$$
\mathbb{P}(\exists n \ge n_0\,:\,g_n \le \log(\delta^{-1})\,|\, Y\left(n_0\right) = X_0 ) = 1\,,
$$ 
and hence almost surely for infinitely many $n$ we have 
$$
\mathrm{d}^{\circ}_{\mathcal{B}}(Y(n)) \ge \delta \sqrt{F(Y(n))}\,.
$$
\end{lemma}
\begin{proof}
 We may assume without loss of generality that $X_0$ has at most one point in $\partial\mathcal{B}$. Indeed,  a.s.~all new points are not in $\partial\mathcal{B}$, and a.s.~$D(Y(n)) \to 0$ as $n \to \infty$, so there exists an $n \ge n_0$ such that $D(Y(n)) < d(Y(n_0))$  and therefore $Y(n) \cap Y(n_0)$ contains at most one point; now apply the lemma with $X_0 = Y(n)$. It follows from these assumptions that a.s.~$g_n \in \mathbb{R}$ for every $n \ge 0$. Likewise we may assume that $\sqrt{\frac{M-1}{M^2}}\sqrt{F(X_0)} < r_0$, where $r_0$ and $c$ are the constants used to specify the uniform geometry of~$\mathcal{B}$, so that by Corollary~\ref{cor: F decreases} and  inequality~\eqref{useful} we have $A(Y(n)) < r_0$ for all $n \ge 0$. 

We will use Lemma~\ref{L: easy isoperimetric} to get an upper bound for the expectation of the positive part of the increment of~$-\log \mathrm{d}^{\circ}_{\mathcal{B}}(Y(\cdot))$.
\begin{equation}\label{E: bounding the positive increments of g} \mathbb{E}\left(\left.\log_{+}\frac{ 
\mathrm{d}^{\circ}_{\mathcal{B}}(Y(n))}{
\mathrm{d}^{\circ}_{\mathcal{B}}(Y(n+1))}\,\right|\, Y(n) = X\right)\,=\,  \frac{\bigintss \mathbf{1}(y \in \mathrm{Keep}(X;\mathcal{B}))\, \log_+\frac{\textstyle{\rule[-4pt]{0pt}{10pt}
\mathrm{d}^{\circ}_{\mathcal{B}}(X) }}{\textstyle{ \rule{0pt}{10pt}\mathrm{dist}\left(y,\mathcal{B}^c\right)}}\,  d\lambda(y)}{\lambda(\mathrm{Keep}(X;\mathcal{B}))}\end{equation}
Let 
$$
r_1 :=\mathrm{d}^{\circ}_{\mathcal{B}}(X),\qquad R := \frac{M+1}{M-1}A(X)\,.
$$
By~\eqref{useful} we have $R \le \frac{M+1}{M\sqrt{M-1}}\sqrt{F(X)}$ and by Lemma~\ref{L: ball containing Keep} we have $\mathrm{Keep}(X;\mathcal{B}) \subseteq B(\mu(X),R) \cap \mathcal{B}$. 
The numerator of the RHS of~\eqref{E: bounding the positive increments of g} may be rewritten as
\begin{align*}
\int \int \mathbf{1}(y \in \mathrm{Keep}(X;\mathcal{B})) &\cdot \mathbf{1}\left(\mathrm{dist}\left(y, \mathcal{B}^c\right) \le r \le
\mathrm{d}^{\circ}_{\mathcal{B}}(X)\right)
\,\frac{dr}{r} \,d\lambda(y) \\ 
& =   \int_0^{r_1} \frac{1}{r} \lambda(\mathrm{Keep}(X;\mathcal{B}) \cap N_r(\mathcal{B}^c))\,dr \\ 
&\le  \int _0^{r_1} \frac{1}{r} \lambda(B(\mu(X),R) \cap \mathcal{B} \cap N_r(\mathcal{B}^c))\,dr \\ 
&\le  
\int _0^{r_1} 2 \sqrt{d} V(d-1) R^{d-1} \,dr\, = 2r_1 d\sqrt{d}\,V(d-1) R^{d-1}\,,
\end{align*}
where the first equality above uses Tonelli's theorem, the next inequality uses Lemma~\ref{L: ball containing Keep}, and the final inequality uses Lemma~\ref{L: easy isoperimetric}.
On the other hand, by our assumption that $A(X) < r_0$, we have  the uniform geometry bound $$\lambda(\mathrm{Keep}(X;\mathcal{B})) \ge c V(d) A(X)^d\,,$$ where $c$ is the constant from~\eqref{E: c}.
Therefore the RHS of~\eqref{E: bounding the positive increments of g} is bounded above by
\begin{eqnarray*}
\frac{2\,\mathrm{d}^{\circ}_{\mathcal{B}}(X)\,  d^{3/2}\,V(d-1) \left(\frac{M+1}{M-1}A(X)\right)^{d-1}}{c V(d) A(X)^d} =  c_1 \frac{
\mathrm{d}^{\circ}_{\mathcal{B}}(X)
}{A(X)} \le  c_1M \frac{\mathrm{d}^{\circ}_{\mathcal{B}}(X)}{\sqrt{F(X)}} \end{eqnarray*} 
where $c_1$ is a constant depending on $d$, $M$ and $c$.  
Thus we have shown that
\begin{equation}\label{E: expected increment of a log-distance} 
\mathbb{E}\left(\left.\log_{+}
\frac{\mathrm{d}^{\circ}_{\mathcal{B}}(Y(n))}{\mathrm{d}^{\circ}_{\mathcal{B}}(Y(n+1))}
\,\right|\, Y(n) = X\right)\,\le\, c_1M \frac{\mathrm{d}^{\circ}_{\mathcal{B}}(X)}{\sqrt{F(X)}}\,. \end{equation}

By almost the same argument, 
$$
%\label{E: small prob of getting nearer to the boundary}
\mathbb{P}(\mathrm{d}^{\circ}_{\mathcal{B}}(Y(n+1))
< \mathrm{d}^{\circ}_{\mathcal{B}}(Y(n))  \, |\, Y(n) = X) 
%\\ 
%\nonumber
%&
\le  \frac{\lambda\left(B(\mu(X),R) \cap \mathcal{B} \cap N_{r_1}\left(\mathcal{B}^c\right) \right)}{\lambda(\mathrm{Keep}(X;\mathcal{B}))} \le c_1M \frac{\mathrm{d}^{\circ}_{\mathcal{B}}(X)}{\sqrt{F(X)}}\,.
$$
Define 
$$
\delta = \frac{1}{8c_1 M^2 4^d},\qquad L = \log(\delta^{-1}) + \frac{1}{2}\log\left(1-\frac1{4M}\right)\,,
$$
(note that $L < \log(\delta^{-1})$) and the stopping time 
$$
\tau_g:= \inf\{n \ge n_0\,:\, g_n \le \log(\delta^{-1})\}\,.
$$  
We aim to prove that $L \vee g_{n\wedge \tau_g}$ is a supermartingale. Noting that $\log(1-1/(4M)) < 0$ and $\log F(Y(n+1)) - \log F(Y(n)) \le 0$, inequality~\eqref{E: probability of definite decrease in F} in Corollary~\ref{cor: significant decrease of F} implies
\begin{align*}
\mathbb{E}&(\max(\log(F(Y(n+1))) - \log F(Y(n)),\,2(L -\log(\delta^{-1}))) \mid \mathcal{F}_n)  \\ 
&= \mathbb{E}\left(\max\left(\log(F(Y(n+1))) - \log F(Y(n)),\,\log\left(1-\frac1{4M}\right)\right)\mid \mathcal{F}_n\right)
\\ 
& 
\le
\mathbb{E}\left(
\mathbf{1}\left(\left.\log(F(Y(n+1))) - \log F(Y(n))\le \log\left(1-\frac1{4M}\right)\right)
\cdot \log\left(1-\frac1{4M}\right)\right| \mathcal{F}_n\right)
\\ 
& = \log\left(1-\frac1{4M}\right)\mathbb{P}\left(\left.\! F(Y(n+1)) - F(Y(n)) \le -\frac{1}{4M}F(Y(n))\,\right|\,\mathcal{F}_n\right)  \le  -\frac{1}{4M}4^{-d}\,. 
\end{align*}
Combining this with~\eqref{E: expected increment of a log-distance}, we find that on the event $g_n \ge \log(\delta^{-1})$, we have 
\begin{align*}
\mathbb{E}((L \vee &g_{n+1}) - (L \vee g_n) \,|\mathcal{F}_n)  =  \mathbb{E}\left(\left(L \vee g_{n+1}\right) - g_n \,|\, \mathcal{F}_n\right) \\&=\, \mathbb{E}\left(\left(L - g_n) \vee (g_{n+1}-g_n)\right) \,|\, \mathcal{F}_n\right) \\&\le\, \mathbb{E}\left((L + \log \delta) \vee (g_{n+1} - g_n) \,|\, \mathcal{F}_n\right)\\ &=\,\mathbb{E}\left(\frac{1}{2}\log\left(1-\frac{1}{4M}\right) \vee (g_{n+1} - g_n) \,|\, \mathcal{F}_n\right)
\\ &\le\, \mathbb{E}\left(\left.\frac{1}{2}\log\left(1-\frac{1}{4M}\right) \vee \left(\frac{1}{2}\log\frac{F(Y(n+1))}{F(Y(n))} + \log_{+}\frac{\mathrm{d}^{\circ}_{\mathcal{B}}(Y(n))}{\mathrm{d}^{\circ}_{\mathcal{B}}(Y(n+1))}\right)  \,\right|\, \mathcal{F}_n\right)\\  
&\le
\,  \mathbb{E}\left(\left.\left( \frac{1}{2}\log\left(1-\frac{1}{4M}\right) \vee \frac{1}{2}\log\frac{F(Y(n+1))}{F(Y(n))}  \right) + \log_{+}\frac{\mathrm{d}^{\circ}_{\mathcal{B}}(Y(n))}{\mathrm{d}^{\circ}_{\mathcal{B}}(Y(n+1))} \,\right|\,\mathcal{F}_n\right)\\  
&\le\,  -\frac{1}{8M}4^{-d} + c_1 M \frac{\mathrm{d}^{\circ}_{\mathcal{B}}(Y(n))}{\sqrt{F(Y(n))}} 
= -\frac{4^{-d}}{8M} + c_1 M e^{-g_n}\,\le\, -\frac{4^{-d}}{8M} + c_1 M \delta\, = 0\,.
\end{align*}
(Note that we used the fact that $\max(a+b,c)\le \max(a,b)+c$ whenever $c\ge 0$ when we went from the fifth to the sixth line.)
We therefore have
\begin{equation}\label{E: g supermartingale}
\mathbb{E}((L \vee g_{(n+1) \wedge \tau_g}) - (L \vee g_{n \wedge \tau_g}) \,|\,\mathcal{F}_n)  \le 0\,.
\end{equation}
Using the first part of Corollary~\ref{cor: significant decrease of F}, we find
$$
\mathbb{P}(g_{n+1} - g_n \le \log(1-1/(4M)) \,|\,\mathcal{F}_n ) \ge 4^{-d} - c_1 M e^{-g_n}\,, 
$$ 
so under the same condition that $g_n \ge \log(\delta^{-1})$ we also have
\begin{equation}\label{E: positive probability that g jumps down a lot}
\mathbb{P}\left((L \vee g_{n+1}) - (L \vee g_n) \le \log\left(1-1/({4M})\right)\,|\,\mathcal{F}_n\right) \ge 4^{-d}\left(1- \frac1{8M}\right) > 0\,.
\end{equation}
Inequality~\eqref{E: g supermartingale} shows that the random sequence $L \vee g_{n \wedge \tau_g}$ is a supermartingale that is bounded below by $L$, so it almost surely converges, but inequality~\eqref{E: positive probability that g jumps down a lot} shows that it can only converge to a value less than or equal to~$\log(\delta^{-1})$. Therefore $\tau_g < \infty$ a.s., and indeed a.s.~$g_n \le \log(\delta^{-1})$ for infinitely many $n$.
\end{proof}

Recall the definition of $\mathrm{d}_{\mathcal{B}}(\cdot)$ from before Lemma~\ref{L: second supermartingale}.
\begin{lemma}\label{L: relative distance to boundary can grow enough} Let\footnote{Think of $\delta$ small and $\Delta$ large.} $\delta > 0$ and $\Delta > 0$. Then there exist $n_1 \in \mathbb{N}$ and $\epsilon > 0$ (both depending on~$d, M, \delta$ and $\Delta$) such that for any $n \in \mathbb{N}$, 
$$ 
\mathbb{P}\left(\frac{\mathrm{d}_{\mathcal{B}}(Y(n+n_1))}{\sqrt{F(Y(n+n_1))}} > \Delta \, \left| \, \frac{\mathrm{d}_{\mathcal{B}}(Y(n))}{ \sqrt{F(Y(n))}} > \delta\right.\right) > \epsilon\,. 
$$
\end{lemma}
\begin{proof}
Before we proceed with the formal proof, let us explain the idea behind it. The main ingredient is that there is a positive probability $\epsilon$ that for some number $n_1$ of consecutive steps the new point arrives very close to the centre of mass of the existing configuration, relative to its diameter. When this occurs, the moment of inertia decreases by a definite multiplicative factor at each step. As a result, after $n_1$ steps the square root of the moment of inertia will have decreased by at least a factor $\delta/(2\Delta)$. At the same time, the convex hull of the configuration stays inside a bounding region that grows from step to step, by a geometrically decreasing sequence of increments. It can be arranged that the bounding region never gets closer to $\mathcal{B}^c$ than $\frac{1}{2}\mathrm{d}_{\mathcal{B}}(Y(n))$, and hence $\mathrm{d}_{\mathcal{B}}(Y(n+n_1)) \ge \frac{1}{2} \mathrm{d}_{\mathcal{B}}(Y(n))$. 

We will combine ideas from the proofs of Corollary~\ref{cor: significant decrease of F} and Lemma~\ref{L: relative tightness of xi}.  Define 
$$
\alpha := \min\left(\frac{M-1}{2M(M+1)}\,,\, \frac{\delta}{48(M+1)}\right)
,\quad
\g:=\sqrt{1-\frac{1}{12(M+1)}}\in(0,1)
$$ 
and 
$$
n_1 := \left\lceil \frac{\log\left(\frac{\delta}{2\Delta}\right)}{\log\g}\right\rceil\,
.
$$ 
For each $i \in \{1, \dots, n_1\}$, define the event $E_i$ by $$ E_i:= \left\{  y_{n+i} \in B\left(\mu(Y(n+i-1)),\,\alpha \sqrt{F(Y(n+i-1))}\right) \right\}\,.$$
We have
$$
\alpha \sqrt{F(Y(n+i-1))} \ge \alpha\frac{M}{\sqrt{M-1}}A(Y(n+i-1)),$$ 
so $$\mathbb{P}(E_i \,|\,Y(n+i-1)) \ge \mathbb{P}\left(y_{n+i} \in  \left. B\left(\mu(Y(n+i-1)), \alpha\frac{M A(Y(n+i-1))}{\sqrt{M-1}} \right)\,\right|\,Y(n+i-1)\right)\,.$$
We have $\frac{\alpha M}{\sqrt{M-1}} < 1$ because $\alpha \le \frac{M-1}{2M(M+1)}$. Hence, using Lemmas~\ref{L: ball containing Keep},~\ref{L: ball inside Keep} and~\ref{L: volume comparison} (as in the proof of Corollary~\ref{cor: significant decrease of F}), we have
$$
\mathbb{P}\left(\left. E_i\,\right|\, Y(n+i-1)\right) \ge \left(\frac{\alpha M\sqrt{M-1}}{M+1}\right)^{d}\,.$$
Therefore 
$$
\mathbb{P}\left(\left. \bigcap_{i=1}^{n_1} E_i \,\right|\, Y(n)\right) \ge \left(\frac{\alpha M\sqrt{M-1}}{M+1}\right)^{d n_1} =: \epsilon\,.$$

Now suppose that $\mathrm{d}_{\mathcal{B}} (Y(n)) \ge \delta \sqrt{F(Y(n))}$ and suppose that the event $\bigcap_{i=1}^{n_1} E_i$ occurs. Our goal is to show that  $\mathrm{d}_{\mathcal{B}} (Y(n+n_1)) \ge \Delta \sqrt{F(Y(n+n_1))}$.

Using the condition $\alpha \le \frac{M-1}{2M(M+1)}$ together with $\sqrt{F(Y(n+i-1))} \le M A(Y(n+i-1))$, see~\eqref{useful}, we obtain for each $i = 1, \dots, n_1$ that 
$$
\alpha \sqrt{F(Y(n+i-1)} \,<\, \frac{M-1}{2(M+1)}A(Y(n+i-1))\,.
$$ 
Then by inequality~\eqref{eq: FF12} from the proof of Corollary~\ref{cor: significant decrease of F}, we have 
$$
F(Y(n+i)) < \g^2\, F(Y(n+i-1))  \,.
$$ 
By induction we find that for $i=1, \dots, n_1$, 
$$
F(Y(n+i)) < \g^{2i}\,F(Y(n))
$$ and in particular 
\begin{equation}\label{E: purple} 
\sqrt{F(Y(n+n_1))} \,\le\, \frac{\delta}{2\Delta} \sqrt{F(Y(n))}\,.
\end{equation}
For every $i = 1, \dots, n_1$, since $\mu(Y(n+i-1)) \in \mathrm{Conv}(Y(n+i-1))$, 
$$
y_{n+i}  \in N_{\alpha\sqrt{F(Y(n+i-1))}}(\mathrm{Conv}(Y(n+i-1)))\,,
$$ 
so
$$ 
\mathrm{Conv}(Y(n+i)) 
\subseteq N_{\alpha \g^{i-1}\sqrt{F(Y(n))}}(\mathrm{Conv}(Y(n+i-1)))\,.
$$ 
Since $$ \sum_{i=1}^{n_1}\alpha \g^{i-1}\, \sqrt{F(Y(n))}  < \beta \sqrt{F(Y(n))}\,,$$ where $$ \beta:= \frac{\alpha}{1-\g}$$
it follows that
$$ 
\mathrm{Conv}(Y(n+n_1)) \subseteq N_{\beta\sqrt{F(Y(n))}}(\mathrm{Conv}(Y(n)))\,.
$$
Moreover, 
$$
\beta \,<\, 24(M+1)\,\alpha\, < \,\frac{\delta}{2}\,,
$$ 
so
\begin{eqnarray}\notag 
\mathrm{d}_{\mathcal{B}} 
(Y(n+n_1)) & \ge & \mathrm{d}_{\mathcal{B}} (\mathrm{Conv}(Y(n))) - \beta\sqrt{F(Y(n))}   =  
\mathrm{d}_{\mathcal{B}} (Y(n)) - \beta\sqrt{F(Y(n))} 
\\
&> &  (\delta - \delta/2)\sqrt{F(Y(n))}\,.\label{E: mauve}
\end{eqnarray}
Taking square roots of both sides of~\eqref{E: purple} and dividing by~\eqref{E: mauve} we obtain
$$ 
\frac{\mathrm{d}_{\mathcal{B}}(Y(n+n_1))}{\sqrt{F(Y(n+n_1))}} > \Delta \,,$$ as required. 
\end{proof}
\begin{lemma}\label{L: get off the boundary} Let $\delta > 0$. Then there exist $n_1' \in \mathbb{N}$ and $\epsilon' > 0$ (both depending on $d, M$, and  $\delta$) such that for any set $X$ of~$M$ points in $\mathcal{B}$ such that
$ |X\cap \partial\mathcal{B}| = 1 \textup{ and } \mathrm{d}_{\mathcal{B}}^\circ(X) > \delta  \sqrt{F(X)}$, and any $n \in \mathbb{N}$, we have
$$ 
\mathbb{P}\left( Y(n+n_1') \subset \mathcal{B}^\circ \, \left| Y(n)=X\right.\right) > \epsilon'\,. 
$$
\end{lemma}
\begin{proof}
The proof is very similar to the proof of Lemma~\ref{L: relative distance to boundary can grow enough}. We take $\Delta = \frac{1}{\sqrt{M-1}}$ and then make the same choices of constants in terms of $M,d, \delta$ and $\Delta$ as we did there:
$$
\alpha := \min\left(\frac{M-1}{2M(M+1)}\,,\, \frac{\delta}{48(M+1)}\right)
,\quad
\g:=\sqrt{1-\frac{1}{12(M+1)}}\,,
$$  
$$
n_1' := \left\lceil \frac{\log\left(\frac{\delta}{2\Delta}\right)}{\log\g}\right\rceil\,
, \quad
\epsilon' = \left(\frac{\alpha M\sqrt{M-1}}{M+1}\right)^{d n_1'}\,.
$$
Let $y_b$ be the unique point in $Y(n) \cap \partial\mathcal{B}$. This time, the good sequence of events whose probability is at least $\epsilon$ is  
$$ E_i' := \left\{ y_{n+i} \in B\left(\mu(Y(n+1-i)\setminus \{y_b\}), \alpha\sqrt{F(Y(n+i-1))}\right) \right\}\,.$$
Instead of applying Lemma~\ref{L: ball inside Keep}, we apply Lemma~\ref{L: other balls inside Keep} to see that  
$$
\mathcal{B} \cap B\left(\mu\left(Y(n+i-1)\setminus \{y_b\}\right), \alpha \frac{MA(Y(n+i-1))}{\sqrt{M-1}} \right) \subseteq \mathrm{Keep}(Y(n+i-1); \mathcal{B})\,.$$ This requires us to check that $\frac{\alpha M}{\sqrt{M-1}} < \frac{M}{M+1}$, which follows from $\alpha \le \frac{M-1}{2M(M+1)}$. It follows as before that on the event $Y(n) \cap \partial\mathcal{B} = \{y_b\}$ we have  
$$ 
\mathbb{P}\left(\left.\bigcap_{i=1}^{n_1'} E_i' \,\right|\, Y(n)\right) \ge \epsilon'.
$$
It also follows as before that on the events $
\mathrm{d}_\mathcal{B}(Y(n) \setminus \{y_b\}) \ge \delta$ and $\bigcap_{i=1}^{n-1} E_i'$ we have

$$ 
\sqrt{F(Y(n+n_1'))} \le \frac{\delta}{2\Delta}\sqrt{F(Y(n))}\,.
$$

Throughout the rest of the proof, we replace each set of the form $\mathrm{Conv}(Y(k))$ by $\mathrm{Conv}(Y(k) \setminus \{y_b\})$. Note that $\mu(Y(k) \setminus \{y_b\}) \in \mathrm{Conv}(Y(k) \setminus \{y_b\})$ so that on the events $E_i'$ we obtain
$$\mathrm{Conv}(Y(n+n_1')\setminus\{y_b\}) \subseteq N_{\beta\sqrt{F(Y(n))}}\mathrm{Conv}(Y(n) \setminus \{y_b\}) $$ and it follows that
$$
\mathrm{d}_{\mathcal{B}}(Y(n+n_1') \setminus \{y_b\}) \ge \frac{\delta}{2}\sqrt{F(Y(n))}\,,
$$ 
as before, and hence
\begin{equation}\label{E: forced off the boundary} \frac{
\mathrm{d}_{\mathcal{B}}^\circ( Y(n+n_1'))}{\sqrt{F(Y(n+n_1'))}} > \Delta = \frac{1}{\sqrt{M-1}}.
\end{equation}

On the other hand, if $y_b \in Y(n+n_1')$ then 
\begin{eqnarray*}F(Y(n+n_1'))
& \ge & \sum_{x \in Y(n+n_1')\setminus\{y_b\}} \Vert y_b - x\|^2\\ & = & (M-1)\left\Vert y_b - \mu(Y(n+n_1')\setminus\{y_b\})\right\Vert^2 \; + \hspace{-1em} \sum_{x \in Y(n+n_1')\setminus\{y_b\}}\hspace{-1em} \left\Vert x - \mu(Y(n+n_1')\!\setminus\!\{y_b\})\right\Vert^2\\ & \ge &  (M-1)\Vert y_b - \mu(Y(n+n_1')\!\setminus\!\{y_b\})\Vert^2 \\ & \ge & (M-1)\,\left[\mathrm{dist}(\mu(Y(n+n_1')\setminus\{y_b\}), \partial\mathcal{B})\right]^2 
\\ & \ge & 
(M-1)\,\left[\mathrm{d}_{\mathcal{B}}\left(\mathrm{Conv}(Y(n+n_1') \setminus \{y_b\})\right)\right]^2 
\\ 
& = & (M-1)\,\left[
\mathrm{d}_{\mathcal{B}}^\circ
(Y(n+n_1'))\right]^2\,,
\end{eqnarray*}
contrary to~\eqref{E: forced off the boundary}. We have shown that on the event $
\mathrm{d}_{\mathcal{B}}
(Y(n)\! \setminus\! \{y_b\}) \ge \delta$, the event $\bigcap_{i=1}^{n-1} E_i'$ occurs with probability at least $\epsilon'$, and when it does, we have $y_b \not\in Y(n+n_1')$, as required.
\end{proof}

\begin{coro}\label{cor: get away from the boundary}
Suppose $\mathcal{B}$ is a convex body in $\mathbb{R}^d$ with uniform geometry and let $Y(n)_{n \ge 0}$ be a $\mathcal{B}$-valued Jante's law process. There exists $\delta > 0$ (depending on $d$ and $\mathcal{B}$) such that for any set $X$ of $M$ distinct points of $\mathcal{B}$, conditional on $Y(n_0) = X$, a.s.~there exists $n  \ge n_0$ such that $
\mathrm{d}_{\mathcal{B}}
(Y(n)) > \delta \sqrt{F(Y(n))}$.
\end{coro}
\begin{proof}
Consider a $\mathcal{B}$-valued Jante's law process started at $Y(n_0)  = X$. Define the stopping time $$\tau_b = \inf\{n \ge n_0 \,:\, Y(n) \cap \partial\mathcal{B} = \emptyset\}\,.$$ 
Note that if $n \ge \tau_b$ then $g_n = \frac{1}{2} \log F(Y(n))  - \log \mathrm{d}_{\mathcal{B}}(Y(n))$. Let $\delta$ be as in Lemma~\ref{L: second supermartingale}, so that a.s.~infinitely often $g_n < \log(\delta^{-1})$. Note that if $n \ge \tau_b$ and $g_n < \log(\delta^{-1})$ then $
\mathrm{d}_{\mathcal{B}}(Y(n)) > \delta \sqrt{F(Y(n))}$. 
So it suffices to show that a.s.~$\tau_b < \infty$. 

We define a sequence of stopping times $\left(\sigma_i\right)_{ i \ge 0}$. Let 
$$
\sigma_0 = \inf\{n \ge n_0\,:\, |Y(n) \cap \partial \mathcal{B} | \le 1\}\,.
$$ 
Since $D(Y(n)) \to 0$ a.s., and a.s.~all new points are in $\mathcal{B}^\circ$, we have a.s.~$\sigma_0 < \infty$. Let $n_1$ and $\epsilon$ be as provided by Lemma~\ref{L: get off the boundary}. Now inductively define for $i = 1, 2, 3, \dots$
$$
\sigma_i = \begin{cases} \inf\left\{n > \sigma_{i-1}\,:\,

\mathrm{d}_{\mathcal{B}}^\circ(Y(n))
> \delta \sqrt{F(Y(n))} \right\} & \text{if $i$ is odd,} \\ \sigma_{i-1}+n_1 & \text{if $i$ is even.} \end{cases} $$

By Lemma~\ref{L: second supermartingale}, a.s.~$\sigma_i < \infty$ for all $i \ge 1$. Then by Lemma~\ref{L: get off the boundary}, for each $m \ge 1$, on the event $\left|Y\left(\sigma_{2m-1}\right) \cap \partial\mathcal{B}\right| = 1$ we have $$\mathbb{P}\left(Y\left(\sigma_{2m}\right) \subset \mathcal{B}^\circ\,|\,\mathcal{F}_{\textstyle{\sigma_{2m-1}}}\right)\,\ge \epsilon.$$
Hence 
$$
\mathbb{P}\left(Y(\sigma_{2m}) \not\subset \mathcal{B}^\circ\right) < (1-\epsilon)^m\,,
$$ so $\tau_b < \infty $ a.s., as required.
\end{proof}

\begin{coro}\label{cor: get far away from the boundary}
Suppose $\mathcal{B}$ is a convex body in $\mathbb{R}^d$ with uniform geometry and let $Y(n)_{n \ge 0}$ be a $\mathcal{B}$-valued Jante's law process. Let $\Delta > 0$ be given. Then for any set $X$ of $M$ distinct points of~$\mathcal{B}$, conditional on $Y(n_0) = X$, a.s.~there exists $n  \ge n_0$ such that $
\mathrm{d}_{\mathcal{B}}(Y(n))
> \Delta \sqrt{F(Y(n))}$.
\end{coro}
\begin{proof}
We combine Lemma~\ref{L: relative distance to boundary can grow enough} with Corollary~\ref{cor: get away from the boundary}. Consider a $\mathcal{B}$-valued Jante's law process started at $Y(n_0)  = X$. Now let $n_1$ and $\epsilon$ be as given by Lemma~\ref{L: relative distance to boundary can grow enough}, applied to the $\Delta$ specified in the statement of Corollary~\ref{cor: get far away from the boundary}. Define a sequence of stopping times $\left(\kappa_i\right)_{i \ge 0}$ as follows. Set $\kappa_0 = n_0$, and then inductively for $i = 1,2,3, \dots$
$$
\kappa_i = \begin{cases} \inf\left\{n > \kappa_{i-1}\,:\,
\mathrm{d}_{\mathcal{B}}(Y(n)) > \delta \sqrt{F(Y(n))} \right\} & \text{if $i$ is odd,} \\ \kappa_{i-1}+n_1 & \text{if $i$ is even.} \end{cases} 
$$
By Corollary~\ref{cor: get away from the boundary}, a.s.~$\kappa_i < \infty$ for all $i$. By Lemma~\ref{L: relative distance to boundary can grow enough}, for each $m \ge 1$, 
$$
\mathbb{P}\left(\left.\mathrm{d}_{\mathcal{B}}
(Y(\kappa_{2m})) > \Delta \sqrt{F(Y(\kappa_{2m}))}\,\right|\, \mathcal{F}_{\textstyle{\kappa_{2m-2}}}
\right) > \epsilon\,,
$$ 
so 
$$\mathbb{P}\left(\mathrm{d}_{\mathcal{B}}
(Y(\kappa_{2m}))  \le \Delta \sqrt{F(Y(\kappa_{2m}))}
\text{ for all } i < 2m
\right) < (1-\epsilon)^m\,.
$$
By letting $m \to \infty$ and using the continuity of probability, as the RHS converges to zero, we arrive to the required conclusion.
\end{proof}

We can now give the proof of Theorem~\ref{T: original}.
\begin{proof}
In Lemma~\ref{L: relative tightness of xi}, we take $\epsilon = \frac{1}{2}$.   Define 
$$ \Delta = \frac{2}{M\sqrt{M-1}}\left(n_0(1/2) + \frac{1}{1 - \gamma^{1/4}}\right)  + \frac{M+1}{M\sqrt{M-1}}\,.$$
On the event that $\mathrm{d}_{\mathcal{B}}(Y(n))
 > \Delta \sqrt{F(Y(n))}$, we have 
$$ 
\mathbb{P}\left(\left.
\Vert \mu(Y(n+j)) - \mu(Y(n))\Vert + \frac{M+1}{M\sqrt{M-1}}\sqrt{F(Y(n))} \le \Delta\sqrt{F(Y(n))}
\text{ for all }j \ge 0
\,\right|\,\mathcal{F}_n\right) \ge \frac{1}{2}
$$ 
and hence, using~\eqref{useful} to compare $A(Y(n+j))$ with $\sqrt{F(Y(n+j))}$,
$$ 
\mathbb{P}\left(\left.\Vert \mu(Y(n+j)) - \mu(Y(n))\Vert + \frac{M+1}{M-1}A(Y(n)) \le \Delta\sqrt{F(Y(n))} \text{ for all }j \ge 0\,\right|\,\mathcal{F}_n\right) \ge \frac{1}{2}\,.
$$
Now by Lemma~\ref{L: ball containing Keep},
$$ 
\mathbb{P}\left( \mathrm{Keep}\left(Y(n+j); \mathbb{R}^d\right) = \mathrm{Keep}(Y(n+j); \mathcal{B}) \text{ for all }j \ge 0\,|\,\mathcal{F}_n\right) \ge \frac{1}{2}\,.
$$

Define an increasing sequence of stopping times $\left(\theta_i\right)_{i \ge 0}$ as follows. Set $\theta_0 = 0$ and then inductively for $i =1, 2, 3, \dots$
$$
\theta_i = \begin{cases} \inf\left\{n > \theta_{i-1}\,:\,\mathrm{d}_{\mathcal{B}}(Y(n)) > \Delta \sqrt{F(Y(n))} \right\} & \text{if $i$ is odd,} \\ \inf\left\{n > \theta_{i-1}\,:\, \mathrm{Keep}\left(Y(n); \mathbb{R}^d\right) \ne \mathrm{Keep}\left(Y(n); \mathcal{B}\right)\right\} & \text{if $i$ is even.} \end{cases} $$
(If for some $i$ we have $\theta_i = \infty$ then for all $j > i$ we also have $\theta_j = \infty$.)
For each $m \ge 1$, on the event that $\theta_{2m-2} < \infty$, firstly $\theta_{2m-1} < \infty$ a.s., by Corollary~\ref{cor: get far away from the boundary}, and then $\mathbb{P}\left(\theta_{2m} = \infty \,\left|\,\mathcal{F}_{\textstyle{\theta_{2m-1}}}\right.\right) \ge \frac{1}{2}$. Hence a.s.~there is a finite $m$ such that $\theta_{2m-1} < \infty$ but $\theta_{2m} = \infty$.

We have shown that the distribution of $\xi$ is a mixture of distributions, each of which is (for some $m$, $n$, and $X$) the conditional distribution of $\xi$ given $\theta_{2m-1} = n$,  $Y(n) = X$ and $\theta_{2m} = \infty$. Each such distribution is equal to the conditional distribution of $z_\infty$ for the $\mathbb{R}^d$-valued Jante's law process $Z(\cdot)$ started at $Z(0) = X$, conditioned on the positive probability event that for all $j \ge 0$, $\mathrm{Keep}\left(Z(j); \mathbb{R}^d\right) = \mathrm{Keep}\left(Z(j); \mathcal{B}\right)$.  By Theorem~\ref{thm: scale-free limit is continuous}, this distribution is absolutely continuous. Thus the distribution of $\xi$ is a mixture of absolutely continuous distributions and is therefore absolutely continuous.
\end{proof}

\section*{Acknowledgment}
The research of S.V.\ is partially supported by the Swedish Science Foundation grant VR 2019-04173 and the Crafoord Foundation grant no.~20190667. S.V.\ would like to acknowledge the hospitality of the University of Bristol during his visits to Bristol. The research of E.C.\ is supported by the Heilbronn Institute for Mathematical Research. We would like to thank John Mackay for pointing us to the notion of uniform geometry in \cite{Leonardi} and Rami Atar for pointing out the relevant Brownian bees model. Finally, we would like to thank the anonymous referee for many useful comments and recommendations, in particular the suggestion to explain what can be done in the non-convex case.

%\appendix
\section*{Appendix 1: Absolute continuity of probability distributions}

Here are the basic facts which we use about absolutely continuous $\mathbb{R}^d$-valued random variables.
\begin{enumerate}
\item {\bf Definitions.} Let $\mathcal{B}(\mathbb{R}^d)$ be the Borel $\sigma$-algebra of Euclidean space $\mathbb{R}^d$. Let $W: M \to \mathbb{R}^d$ be a Borel-measurable random variable on a probability space $(M, \mathfrak{S},\mathbb{P})$, meaning that for every $A \in \mathcal{B}(\mathbb{R}^d)$, $W^{-1}(A) \in \mathfrak{S}$. The distribution of $W$ is the unique Borel probability measure $\nu$ on $\mathcal{B}(\mathbb{R}^d)$ such that for every $A \in \mathcal{B}(\mathbb{R}^d)$ we have $\nu(A) = \mathbb{P}(W \in A)$. (It suffices to check that this holds for every open set $A \subset \mathbb{R}^d$, or for every closed set.) We call both $W$ and $\nu$ \emph{absolutely continuous} when $\nu$ is absolutely continuous with respect to the Lebesgue measure $\lambda$. This means that for every Borel set $A \subset \mathbb{R}^d$, $$(\lambda(A) = 0) \implies (\nu(A) = \mathbb{P}(W \in A) = 0)\,.$$ Equivalently, there exists an element $f \in L^1(\mathbb{R}^d,\lambda)$ of norm $1$ such that for every Borel set $A \subseteq \mathbb{R}^d$ we have 
$$
\mathbb{P}(W \in A) = \nu(A) = \int_A  f(x) \,d\lambda(x).
$$ 
This element $f$ is called the density of $\nu$, or the density of $W$. Every element of $L^1(\mathbb{R}^d, \lambda)$ has a Borel-measurable representative, so we may assume that $f$ is Borel-measurable. The density of an absolutely continuous random variable is not necessarily representable by a continuous function, and if $M$ happens to be equipped with a topology which generates $\mathfrak{S}$, the absolute continuity of $W$ as a random variable is unrelated to the continuity of $W$ as a function.  

\item {\bf Mixtures of absolutely continuous random variables are absolutely continuous.} Suppose that $\pi$ is a Borel probability measure on $L^1(\mathbb{R}^d)$ and $\nu$ is a random absolutely continuous probability measure defined on $(M, \mathfrak{S},\mathbb{P})$ whose density $d\nu/d\lambda$ with respect to Lebesgue measure $\lambda$ is distributed according to $\pi$. Then for any Lebesgue measurable set $A \subseteq \mathbb{R}^d$ we have (by Tonelli's theorem)
\begin{eqnarray*}
\mathbb{E}^\pi(\nu)(A) := \int \nu(A)\,d\pi(\nu) & = &  \int\int \frac{d\nu}{d\lambda}(x) \mathbf{1}_A(x) \,d\lambda(x)\,d\pi(\nu) \\ & = & \int \int \frac{d\nu}{d\lambda}(x) \,d\pi(\nu) \mathbf{1}_A(x)\,d\lambda(x)\,,
\end{eqnarray*}
so the mixture $\mathbb{E}^\pi(\nu)$ is an absolutely continuous probability measure with density $\int \frac{d\nu}{d\lambda}\,d\pi(\nu)$. To see that the integral of this density with respect to $\lambda$ is $1$, take $A = \mathbb{R}^d$ above.  Even if $\nu$ a.s.~has a $C^\infty$ density, it need not be true that $\mathbb{E}^\pi(\nu)$ has a continuous density.  

\item {\bf Conditioning on positive probability events.} Let $\nu$ be a Borel probability measure on $\mathbb{R}^d$ and $W$ an $\mathbb{R}^d$-valued random variable on  $(M, \mathfrak{S}, \mathbb{P})$ with distribution $\nu$. If $E \in \mathfrak{S}$ is an event such that $\mathbb{P}(E) > 0$, then the conditional distribution of $W$ given $E$ is also absolutely continuous.
Indeed, let $\nu_E(A)=\P(W\in A\mid E)$. Then
$$
\nu_E(A)=\frac{\P(E\cap\{W\in A\})}{\P(E)}
\le \frac{\nu(A)}{\P(E)}
$$
so that $\nu(A)=0$ implies $\nu_E(A)=0$, so that $\nu_E\ll \nu\ll \lambda$ where $\lambda$ is the Lebesgue measure.
\end{enumerate}

\section*{Appendix 2: Proof of Lemma~\ref{L: isoperimetric}}
\begin{proof}
 By replacing $\mathcal{B}$ with $\mathcal{B} \cap B(x,R+r+1)$, which does not alter the set whose volume is to be estimated, we may assume that $\mathcal{B}$ is bounded.
 
Next, we reduce to the case where $\mathcal{B}$ is a polyhedron. Let $\epsilon > 0$. We claim that there is a convex polyhedron $P$ (with finitely many facets) such that $\mathcal{B} \subseteq P \subset N_\epsilon(\mathcal{B})$. To see this, let $y$ be any point in the interior of $\mathcal{B}$ and let $\pi_y$ denote the radial projection about $y$ from $\mathbb{R}^d\setminus\{y\}$ onto the unit sphere $\mathbb{S}_y$ that is centred on $y$. For each point $z \in \partial\mathcal{B}$ and each supporting hyperplane~$H$ of~$\mathcal{B}$ at $z$, the radial projection $\pi_y\left(H \cap N_\epsilon(\mathcal{B})\right)$ is a neighbourhood of $\pi_y(z)$ in the sphere $\mathbb{S}_y$. Since $\mathbb{S}_y$ is compact, some finite collection of such projections covers it. The corresponding finite collection of support hyperplanes defines a suitable polyhedron $P$ (by taking the intersection of the half-spaces that they bound which contain $\mathcal{B}$).
 
For each $z \in \mathcal{B}$, there exists a point $b \in \partial\mathcal{B}$ such that $\Vert z - b \Vert = \mathrm{dist}(z,\mathcal{B}^c)$, and a supporting hyperplane $H$ of $\mathcal{B}$ such that $H$ passes through $b$ and if $z \neq B$ then $H$ is orthogonal to $b-z$. Then the translate of $H$ by distance $\epsilon$ in the direction $b-z$ is a supporting hyperplane of $N_\epsilon(\mathcal{B})$ and it follows that $$  \mathrm{dist}\left(z, N_\epsilon\left(\mathcal{B})^c\right)\right)  = \mathrm{dist}\left(z, \mathcal{B}^c \right) + \epsilon\,.$$ Since $P \subseteq N_\epsilon(\mathcal{B})$, we have
 $$\mathrm{dist}\left(z, P^c\right) \le \mathrm{dist}\left(z, \mathcal{B}^c\right)\,.$$
It follows that
 $$
 \mathcal{B} \cap N_r\left(\mathcal{B}^c\right) \subset P \cap N_{r+\epsilon}\left(P^c\right)\,.
 $$
Hence if inequality~\eqref{E: isoperimetric} always holds for a bounded convex polyhedron $P$ with finitely many faces in place of $\mathcal{B}$, (and $r+\epsilon$ in place of $r$), then by taking the limit as $\epsilon \searrow 0$, we obtain~\eqref{E: isoperimetric} for a general bounded convex body $\mathcal{B}$ and hence also for any convex body $\mathcal{B}$.
 
Finally, we prove~\eqref{E: isoperimetric} in the case where $\mathcal{B}$ is a bounded convex polyhedron $P$ which has finitely many faces.  Suppose the facets of $P$ are $F_1, \dots, F_n$. (A facet is a face of co-dimension one.) For each $i = 1, \dots, n$, let $u_i$ be the outward-pointing normal vector of~$F_i$, and let the facet~$F_i$ be contained in the hyperplane $H_i = \{z \in \mathbb{R}^d \,:\, z.u_i = h_i\}$, so that $$P = \{z \in \mathbb{R}^d \,:\, z.u_i < h_i \text{ for $i =1, \dots, n$}\}\,.$$

 Let $\mathrm{Cut}(P)$ be the \emph{cut-set} of $P$, which is the set of points in $P$ that are equidistant from at least two points in $\partial P$. Note that $\lambda(\mathrm{Cut}(P)) = 0$ since $\mathrm{Cut}(P)$ is contained in the union of the finite set of hyperplanes which bisect the dihedral angles of pairs of hyperplanes $H_i, H_j$. Also, $\lambda(\partial{P}) = 0$. So to prove~\eqref{E: isoperimetric}, if suffices to define a volume-preserving injection $\varphi$ from $B(x,R) \cap (P^\circ\setminus\mathrm{Cut}(P)) \cap N_r(P^c)$ to $B(x,R) \setminus B(x,R-r)$.  
 
Partition the open set $(P^\circ\setminus\mathrm{Cut}(P)) \cap N_r(P^c)$ into finitely many pieces $Q_1, \dots, Q_n$, where each piece $Q_i$ consists of the points whose closest facet of $P$ is a given facet $F_i$. $Q_i$ is an open convex polyhedron contained in the $r$-neighbourhood of the facet $F_i$. 

The map $\varphi$ will translate each point of $Q_i \cap B(x,R)$ by a non-negative multiple of $u_i$, where the multiple may vary from point to point.  Specifically, for any $z \in Q_i \cap B(x,R)$, define 
$$
a(z) = \sup(\{t \in \mathbb{R}\,:\, z + t u_i \in Q_i \}) = h_i - z.u_i = \mathrm{dist}(z, P^c)\,,
$$
and
$$
b(z) = \sup(\{t \in \mathbb{R}\,;\, z + t u_i \in B(x,R)\})\,. 
$$
Then define
 $$ \varphi(z) :=  z + \max(0,b(z) - a(z))\,.$$
Note that if $z' = z + t u_i$ is also in $Q_i \cap B(x,R)$ then $b(z') = b(z) - t$ and $a(z') = a(z) - t$, since both $B(x,R)$ and $Q_i$ are convex. Hence $\varphi(z') = \varphi(z) + t u_i$.
Informally, the map $\varphi$ slides each line segment of the form $\{z + t u_i:  t \in \mathbb{R},  z + t u_i \in B(x,R) \cap Q_i\}$ as far as possible in the direction~$u_i$ (parallel to the line segment) such that it remains a subset of $B(x,R)$. Each of these line segments has length at most $r$. It follows that every point in the image of $\varphi$ is within distance $r$ of the boundary of $B(x,R)$. The restriction of $\varphi$ to $Q_i \cap B(x,R)$ is a smooth volume-preserving map, since both $a(z)$ and $b(z)$ depend smoothly on $z$, and in an appropriate coordinate system its Jacobian matrix is upper unitriangular. 

Finally, we must check that $\varphi$ is an injection. From the description in terms of sliding line segments, we see that the restriction of $\varphi$ to each piece $Q_i \cap B(x,R)$ is an injection. Secondly, if $z \in Q_i \cap B(x,R)$, then the unique closest facet of $P$ to $\varphi(z)$ is $F_i$. Therefore the images of distinct pieces under $\varphi$ are disjoint and we are done.

\end{proof}

\section*{Appendix 3: Beyond convex bodies}

For a bounded measurable set $\mathcal{B} \subset \mathbb{R}^d$, we will explain below how some reasonable geometric assumptions on $\mathcal{B}$, which are more general than convexity, are sufficient to define the $\mathcal{B}$-valued Jante process (so that it is almost surely well-defined for all times) and to show that it converges almost surely to a random limit. However, it appears to be more difficult to show that the limit point has an absolutely continuous distribution, even in the nicest possible non-convex setting, where $\mathcal{B}$ is a region bounded by a compact smooth hypersurface in $\mathbb{R}^d$. In this case, the boundary of $\mathcal{B}$ is Lebesgue null, so to prove absolute continuity of the distribution of the limit point one must prove in particular that the limit point almost surely lies in the interior of $\mathcal{B}$.  In Lemma~\ref{L: non-convex continuity} below, which applies to a reasonable class of non-convex bodies, we show that the \emph{only obstacle} is proving that the limit point almost surely lies in the interior of $\mathcal{B}$.

First, the notion of $\mathcal{B}$-valued Jante process may be extended to any subset $\mathcal{B} \subset \mathbb{R}^d$ satisfying the following condition:
\begin{quote}
{\bf{Condition A}}: $\mathcal{B}$ is Lebesgue measurable and non-empty, and $\overline{\mathcal{B}} = \mathrm{support}(\lambda|_\mathcal{B})$, where $\lambda$ denotes Lebesgue measure on $\mathbb{R}^d$. (In particular this implies that $\lambda(\mathcal{B}) > 0$.) \end{quote}
Condition $A$ holds for a non-empty Lebesgue measurable subset $\mathcal{B} \subset \mathbb{R}^d$ if and only if for every $x \in \mathcal{B}$ and $\epsilon > 0$ we have $\lambda(B(x,\epsilon) \cap \mathcal{B}) > 0$.   For example, condition $A$ holds whenever is $\mathcal{B}$ is the closure of its interior. 
\begin{lemma}
Suppose $\mathcal{B}$ satisfies condition A and $X = \{x_1, \dots, x_M\}$ is a set of $M$ distinct points of $\mathcal{B}$. Then the Lebesgue measure of $\mathrm{Keep}(X;\mathcal{B})$ is positive.
\end{lemma}
\begin{proof} Let $x_j$ be a point of $X$ such that $d(x_j,\mu(X)) = A(X)$, i.e. $x_j$ is a point of $X$ that is as far as possible from $\mu(X)$. Let $x_i$ ($i \neq j$) be another point of $X$. Since $x_i$ lies in the closure of the  ball $B(\mu(X), A(X))$, which is internally tangent to the ball $B\left(\frac{\Sigma(X) - x_j}{M-1}, \frac{M}{M-1}A(X)\right)$ at $x_j$, but $x_i \neq x_j$, we see that $x_i$ lies in the (open) ball $B\left(\frac{\Sigma(X) - x_j}{M-1}, \frac{M}{M-1}A(X)\right)$ and hence there is an $\epsilon > 0$ such that $$B(x_i,\epsilon) \subset B\left(\frac{\Sigma(X) - x_j}{M-1}, \frac{M}{M-1}A(X)\right)\,,$$ and hence by Lemma 2 we have $B(x_i,\epsilon) \cap \mathcal{B} \subset \mathrm{Keep}(X;\mathcal{B})$. Condition A implies that $B(x_i, \epsilon) \cap \mathcal{B}$ has positive Lebesgue measure, so $\mathrm{Keep}(X;\mathcal{B})$ has positive Lebesgue measure, as required.
\end{proof}
Under condition $A$, the $\mathcal{B}$-valued Jante process is well-defined starting at any set $Y(0)$ of $M \ge 2$ distinct points of $\mathcal{B}$. Indeed, when $Y(n)$ is any set of $M \ge 2$ distinct points in $\mathcal{B}$, it makes sense to sample a point uniformly from $\mathrm{Keep}(Y(n); \mathcal{B})$ in order to define $Y(n+1)$, because this set has positive but finite Lebesgue measure. Moreover, almost surely $Y(n+1)$ will also consist of $M$ distinct points. 

 To ensure that the $\mathcal{B}$-valued Jante process almost surely converges to a limit point, it is helpful to assume a stronger condition. 
\begin{defn}
Let $\zeta$ be an $\mathbb{R}^d$-valued random variable. Any subset $A \subset \mathrm{support}(\zeta)$ is called \emph{regular} with parameters $\delta_A \in (0,1), \sigma_A, r_A > 0$ if for every $x \in A$ and $0 < r \le r_A$ we have
$$\mathbb{P}(\zeta \in B(x,r\delta_A) \,|\, \zeta \in B(x,r)) \ge \sigma_A\,.$$
$\zeta$ is called \emph{regular} if for every $x \in \mathrm{support}(\zeta)$, there exists a $\gamma_x > 0$  such that the subset $B_\gamma(x) \cap \mathrm{support}(\zeta)$ is regular (with any parameters). 
\end{defn}
\noindent This regularity condition is Assumption 2 in Kennerberg and Volkov \cite{KV1}. Specializing to the case at hand, we introduce
\begin{quote}
{\bf{Condition B:}}  $\mathcal{B}$ is a bounded measurable subset of $\mathbb{R}^d$ such that $\lambda(B) > 0$ and such that the random variable $\zeta$ whose distribution is $\frac{1}{\lambda(\mathcal{B})}\lambda|_\mathcal{B}$ is regular.
\end{quote}
We remark that condition B is closely related to the uniform geometry condition that we defined for the case of convex sets. Under conditions A and B, the $\mathcal{B}$-valued Jante process is a special case of the process studied by Kennerberg and Volkov, and \cite[Theorem 2]{KV1} implies that the $\mathcal{B}$-valued Jante process almost surely has a (random) limit $y_\infty$ in the closure of $\mathcal{B}$. 

To make our coupling approach to generalize to a non-convex set $\mathcal{B}$, in addition to conditions A and B, we need two more conditions:
\begin{quote} {\bf{Condition C:}} $\mathcal{B}$ is the closure of its interior. \end{quote}
\begin{quote} {\bf{Condition D:}} The limit point $y_\infty$ of the $\mathcal{B}$-valued Jante process almost surely lies in the interior of $\mathcal{B}$. \end{quote}

Unfortunately, we do not know whether \emph{in general} conditions A, B and C imply condition~D.  In some special circumstances other than the convex case, it might be possible to prove a substitute for  Lemma~\ref{L: easy isoperimetric} and hence also a substitute for Corollary~\ref{cor: get far away from the boundary}. For example, if $B$ is bounded by a compact smooth hypersurface, then conditions A, B and C hold, and one might also hope for a substitute for Lemma~\ref{L: easy isoperimetric} to hold with uniform constants because the principal curvatures of the bounding hypersurface would be bounded.

\begin{lemma}\label{L: non-convex continuity} Suppose  $\mathcal{B} \subset \mathbb{R}^d$ satisfies conditions A, B, and C. Consider the $\mathcal{B}$-valued Jante process with initial state $Y(0) \subset \mathcal{B}$ consisting of $M \ge 2$ distinct points of $\mathcal{B}$. Then if condition~D also holds, the distribution of $y_\infty$ is absolutely continuous.\end{lemma}
\begin{proof}
Assume conditions A, B, C, and D hold. Take a countable cover of the interior of $\mathcal{B}$ by distinct open balls, $$\mathcal{C} = \{U_i: i \in \mathbb{N}\}\,,$$ where each $U_i \subseteq \mathcal{B}$.  Then define a random element $U_\infty$ of $\mathcal{C}$ by letting $U_{\infty} = U_i$ where $i \in \mathbb{N}$ is minimal such that $y_\infty \in U_i$. By the assumptions, almost surely at least one such index $i$ exists.
 
For any $U \in \mathcal{C}$ and $T \ge 0$, let $K(U,T)$ be the event that for all $t \ge T$, $\mathrm{Keep}(Y(t);\mathcal{B}) \subset U$.  Let~$T_0$ be the earliest time $T$ such that $K(U_\infty, T)$ holds. (Note that $T_0$ is not a stopping time, but it is almost surely finite.) Then the distribution of $y_\infty$ is a mixture of the conditional distributions of $y_\infty$ 
conditioned on $U_{\infty} = U$ and $T_0 = T$, as $(U,T)$ ranges over the countable set
$$
\{(U,T) \in \mathcal{C} \times \mathbb{N}_0\,:\, \mathbb{P}(U_{\infty} = U, T_0= T) > 0\}\,.
$$
It therefore suffices to prove that each of these conditional distributions is absolutely continuous.  So suppose that that $U \in \mathcal{C}$ and $T \ge 0$ satisfy $\mathbb{P}(U_\infty = U, T_0 = T) > 0$. The conditional distribution of $y_\infty$ given $U_{\infty} = U$ and  $T_0 = T$ is obtained from the conditional distribution of~$y_\infty$ given only $K(U,T)$  by further conditioning on the \emph{positive probability} event that $y_\infty$ is not contained in any $U' \in \mathcal{C}$ which precedes $U$ in the enumeration. So it suffices to prove that the conditional distribution of $y_\infty$ given $K(U,T)$ is absolutely continuous. 

The conditional distribution of~$y_\infty$ given $K(U,T)$ is a mixture of the distributions of $y_\infty$ given~$Y(T)$ and~$K(U,T)$, where the mixture is with respect to the distribution of $Y(T)$ conditioned on the positive probability event $K(U,T)$. So it suffices to prove the absolute continuity of the conditional distribution of $y_\infty$ given $Y(T) = X$ and $K(U,T)$, for every configuration $X$ such that $\mathbb{P}(K(U,T) \,|\, Y(T) = X) > 0$.

Using a coupling argument, we will show that this distribution is absolutely continuous with respect to the distribution of the limit point of $U$-valued Jante process started at $X$, which is absolutely continuous by our Theorem~\ref{T: original}. Since the Jante process is homogeneous, we may replace $T$ with $0$ without loss of generality in the rest of the proof.

Let $Y'$ be a $U$-valued Jante process started in configuration $X \subset U$, and let $Y$ be a $\mathcal{B}$-valued Jante process started in configuration $X$. Define the events
$$K(U) : = \{(\forall t \ge 0)(\mathrm{Keep}(Y(t); \mathcal{B}) \subset U)\}$$
and
$$K'(U) := \{(\forall t \ge 0)(\mathrm{Keep}(Y'(t); \mathcal{B}) \subset U)\}\,.$$ We claim that $\mathbb{P}(K'(U)) = \mathbb{P}(K(U))$. Indeed, there is a coupling of $\left(Y(t)\right)_{t \ge 0}$ and $\left(Y'(t)\right)_{t \ge 0}$ such that $Y(t) = Y'(t)$ for all times $t$ up to and including the first time when $\mathrm{Keep}(Y(t);\mathcal{B}) \not\subset U$, or for all $t$ if there is no such time. (After such a time, we may let $Y'$ and $Y$ evolve independently.) Let $y'_\infty$ be the limit point of $Y'(\cdot)$, which exists almost surely. Let $\mu_{U,X}$ be the (defective) positive Borel measure defined by
\begin{eqnarray*}\mu_{U,X}(\cdot) & = & \mathbb{P}(y_\infty \in \cdot \text{ and } K(U) )\\
& = & \mathbb{P}(y'_\infty \in \cdot \text{ and } K'(U) )\,.
\end{eqnarray*}
(The equality here is a consequence of the coupling.) 
Since $\mathbb{P}(K(U)) = \mathbb{P}(K'(U))$, we deduce that the conditional distribution of~$y_\infty$ given~$K(U)$ is equal to the conditional distribution of~$y'_\infty$ given~$K'(U)$. (Indeed, both are obtained by normalizing~$\mu_{U,X}$.) 
Finally, note that $\mu_{U,X}$ is absolutely continuous with respect to the unconditional distribution of $y'_\infty$, with the Radon-Nikodym derivative at $y$ being $\mathbb{P}(K'(U) \,|\, y'_\infty = y)$.  Finally, the unconditional distribution of $y'_\infty$ is absolutely continuous, by Theorem~\ref{T: original}.

Putting this all together we conclude that $y_\infty$ has an absolutely continuous distribution.  
\end{proof}
\end{document}